\documentclass[12pt, reqno]{amsart}
  
\usepackage{ stmaryrd }
\usepackage{cite}

\usepackage{caption}

\usepackage{color}
 \usepackage{amscd}

\usepackage{amsrefs}

 \usepackage{amssymb, amsmath}

 \usepackage[final]{hyperref}

 \input xy
  \xyoption{all}
 
 \input diagxy

\usepackage{tikz, tikz-cd}
\usepackage{tikz-3dplot}

\usetikzlibrary{decorations.markings}

\usetikzlibrary{shapes,snakes}

\usetikzlibrary{arrows.meta}

\tikzstyle{mybox} = [draw=black,   
    rectangle, rounded corners, inner sep=5pt, inner ysep=4pt]

\tikzset{middlearrow/.style={
        decoration={markings,
            mark= at position 0.5 with {\arrow{#1}} ,
        },
        postaction={decorate}
    }
}

\newcommand{\tikzAngleOfLine}{\tikz@AngleOfLine}
  \def\tikz@AngleOfLine(#1)(#2)#3{%
  \pgfmathanglebetweenpoints{%
    \pgfpointanchor{#1}{center}}{%
    \pgfpointanchor{#2}{center}}
  \pgfmathsetmacro{#3}{\pgfmathresult}%
  }

\usepackage{epsf}
\usepackage{epsfig}

 \newtheorem{lem}{Lemma}

    \newtheorem{prop}{Proposition}

  \numberwithin{equation}{section}
   \numberwithin{lem}{subsection}

 \numberwithin{prop}{subsection}
 \numberwithin{theorem}{subsection}

\newcommand{\Beq}{\begin{equation}}
\newcommand{\Eeq}{\end{equation}}
\newcommand{\Beqr}{\begin{eqnarray}}
\newcommand{\Eeqr}{\end{eqnarray}}

\newcommand{\mbm}{{\mathbf  M}}

\newcommand{\mbs}{{\mathbb  S}}

\newcommand{{\mbbx}}{{\mathbb  X}}

\newcommand{{\mbx}}{{\mathbf  X}}

\newcommand{\mbr}{{\mathbb R}}

\newcommand{\mcg}{{\mathcal G }}

\newcommand{\mcp}{{\mathcal P }}

\newcommand{\mck}{{\mathcal K }}

\newcommand{\ovgam}{  {\overline \gamma}}

\newcommand{\ovdel}{  {\overline \delta}}

\newcommand{{\wtlg}}{\widetilde\gamma }

\newcommand{{\wtlG}}{\widetilde\Gamma }

\newcommand{{\wtlv}}{\widetilde v }

\newcommand{\tld}{\tilde\delta }

\newcommand{\tlg}{\tilde\gamma }

\newcommand{{\tlG}}{\tilde\Gamma }

\newcommand{\mbbad}{\mbba^{\bullet\bullet}}

\newcommand{\mbba}{\mathbf {A }}

\newcommand{\mbg}{\mathbf {G} }
\newcommand{\mbgd}{\mbg^{\bullet\bullet}}

\newcommand{\mbp}{\mathbf {P} }

\newcommand{\mbq}{\mathbf {Q} }

\newcommand{\mbpd}{\mbp^{\bullet\bullet}}

\newcommand{\Obj}{{\rm Obj }}
\newcommand{\Mor}{{\rm Mor }}

\newcommand{\mbbs}{\mathbb{S} }

\begin{document}
\title[Pushforwards]{Pushforwards and  Gauge Transformations for Categorical Connections}

\author{Saikat Chatterjee }
\address{School of Mathematics\\
Indian Institute of Science Education and Research - Thiruvananthapuram\\
 Maruthamala PO, Vithura,\\ Thiruvananthapuram 695551, Kerala, India}

\email{saikat.chat01@gmail.com}
\author{Amitabha Lahiri}
\address{Amitabha Lahiri, S.~N.~Bose National Centre for Basic Sciences \\ Block JD,
  Sector III, Salt Lake, Kolkata 700098 \\
  West Bengal, India}
  \email{amitabhalahiri@gmail.com}

\author{Ambar N. Sengupta }
\address{Ambar N. Sengupta, Department of Mathematics\\
 Department of Mathematics\\
  University of Connecticut\\
 Storrs, CT 06269, USA}
\email{ambarnsg@gmail.com}

 

\def\xypic{\hbox{\rm\Xy-pic}}

\begin{abstract}   We construct and study pushforwards of categorical connections on categorical principal bundles. Applying this construction to the case of decorated path spaces in principal bundles, we obtain a transformation of classical connections that combines the traditional gauge transformation with an affine translation. \end{abstract}


 \maketitle
 
 \section{Introduction}

In classical gauge field theories the field of interest is described  mathematically by a Lie-algebra-valued $1$-form   $A$ that is subject to   gauge transformations  $A\mapsto \theta^{-1}A\theta +(d\theta)\theta^{-1}$. Extending this idea,  some physical fields could be described by higher order Lie-algebra-valued forms along with suitable transformation laws for these forms. Such forms can be viewed in terms of parallel transport processes over surfaces and higher dimensional submanifolds of the ambient manifold. Parallel transport over spaces of paths have been of interest in mathematics \cites{Gross1985, Sing1995}, inspired by problems arising from Yang-Mills theory.  Alvarez et al. \cite{Alva1998} studied ``higher'' parallel transport in the context of several physical theories such as BF theories, Chern-Simons theory, and the self-dual Yang-Mills equations.

Categorical bundle theory provides a framework in which such ``higher'' parallel transport processes can be formalized.  There are distinct formalisms for categorical bundle  theory.  In this paper we follow the categorical framework for connections over path spaces developed in  \cites{CLS2geom, CLSpdg}; for ease of reference, section \ref{s:cpb} includes a largely self-contained description of the framework.

The main objectives of this paper are to:
\begin{enumerate}
\item Construct and study the notion of a {\em pushforward}  for   connections on categorical bundles;
\item Use the pushforward to construct an extension of the notion of the traditional gauge transformation to include affine translates $A\mapsto \phi^*A+\Lambda$, where $\phi$ is a traditional gauge transformation and $\Lambda$ is an appropriate type of $1$-form.
\end{enumerate}
A principal $G$-bundle is a smooth submersion $\pi:P\to M$ of manifolds, along with a  smooth free right action $(p,g)\mapsto pg$ of a Lie group $G$ on $P$, preserving the fibers  of the projection $\pi$; there is also a local triviality property (see Kobayashi and Nomizu \cite{KobNomI} for the theory).    A connection $A$ on this bundle is a $1$-form on $P$ with values in the Lie algebra $L(G)$ with certain properties; the geometric significance of $A$ is that it leads to a way of lifting a path $\gamma$ on $M$ to a path ${\tlg}_p$, initiating at $p$, on $P$, with $\pi\circ\tlg_p=\gamma$. There are different counterparts of this theory for the categorical context. In the approach we follow, a categorical principal bundle is given by a functor $\pi_{\mbp}:\mbp\to\mbm$, and there is a an action $\mbp\times\mbg\to\mbm$, where $\mbg$ is a categorical group; we will explain these concepts in section \ref{s:cpb}. A categorical connection is a prescription to lift morphisms $\gamma$  of the base category $\mbm$ to morphisms $\ovgam_p$ in $\mbp$.  A traditional gauge transformation is a smooth mapping $P\to P$ that preserves fibers and the action of $G$, and is specified by a smooth function $\theta:P\to G$ which has an equivariance property. We will show (in section \ref{ss:gtdiff}) that the categorical counterpart of this, a {\em categorical gauge transformation},  is specified by both {\em the function $\theta$ and a $1$-form on $\Mor(\mbp)$ that takes values in the Lie algebra of} a subgroup of $\Mor(\mbg)$.

\subsection{Technical description} In classical bundle theory a connection on a principal bundle can be pushed forward to produce a connection on a different bundle. In more detail, suppose $\pi_P:P\to M$ and $\pi_Q:Q\to M$ are principal $G$- and $K-$bundles, where $G$ and $K$ are Lie groups. Suppose $s:G\to K$ is a Lie group homomorphism and 
 \begin{equation}\label{E:PQmapintro}
\begin{tikzcd}[column sep={3em,between origins}]
P\arrow{rr}{S}\arrow{dr}[swap]{\pi_{P}} &&
Q   \arrow{dl}{\pi_{Q}}   \\ 
& M  
\end{tikzcd}
\end{equation}
a commutative diagram, with $S$ a smooth map that satisfies $S(pg)=S(p)s(g)$ for all $p\in P$ and $g\in G$. Then a connection $A_P$ on $P$ produces a connection $S_*A_P$ on $Q$ essentially by declaring that the $S$ map $A_P$-horizontal paths on $P$ to $S_*A_P$-horizontal paths on $Q$. 

A basic example is the case where $M$ is a Riemannian manifold of dimension $n$, $K$ is the orthogonal group $O(n)$, with $s$ being the inclusion into the general linear group $GL(n)$,    $P$ is the bundle of orthonormal frames for the tangent bundle of $M$, and $S$ is the inclusion map into the bundle  $Q$ of all frames over $M$. 

In section \ref{s:cpb} we give a self-contained description  of the mathematical formalism developed in our earlier works, but in a modified form that brings out some features more clearly. We review the notion of a categorical group $\mbg$, which involves two groups $G$ and $H$ intertwined in a special structure called a Lie crossed module.  We also describe the notions of  categorical path spaces and a categorical principal bundle $\pi_{\mbp}:\mbp\to\mbm$.  Briefly, $\mbm$ is a category whose objects are points of a manifold and whose morphisms correspond to paths on the manifold. A central example of interest for categorical bundles is that of a {\em decorated} bundle $\mbp^{A, {\rm dec}}\to\mbm$, which arises from a classical principal $G$-bundle $\pi:P\to M$, a connection $A$ on this bundle, and an additional new structure group $H$ as mentioned above. Then the objects of  $\mbp^{A, {\rm dec}}\to\mbm$ are just the points of $P$, while morphisms are of the form $(\tlg^A, h)$, where $\tlg^A$ is any $A$-horizontal path on $P$ and $h\in H$ is a decoration of that path.

 In section \ref{s:pfcon} we construct and study the categorical counterpart of the pushforward (\ref{E:PQmapintro}). Briefly, if
 \begin{equation}\label{E:PQmapcatintro}
  \begin{tikzcd}[column sep={3em,between origins}]
\mbp\arrow{rr}{\mbs}\arrow{dr}[swap]{\pi_{\mbp}} &&
\mbq  \arrow{dl}{\pi_{\mbq}}   \\ 
& \mbm 
\end{tikzcd}
\end{equation}
is a categorical counterpart to the diagram (\ref{E:PQmapintro}) then, by pushing forward horizontal lifts, we obtain a categorical connection $\mbs_*\mbba_{\mbp}$ on $\mbq$ from a given categorical connection $\mbba_\mbp$ on $\mbp$.
 
 In our presentation we build up to this general notion of pushforward by  {\em first studying  the examples } of interest in sections \ref{ss:functmbs}  and  \ref{ss:projcon}.

In section \ref{ss:liftcon} we use a process that is a kind of inverse of the pushforward that works only in the context we need. Briefly, if $\mbs:\mbp\to\mbq$ is a functor between categorical principal bundles, preserving all relevant structures,  then a categorical connection   $\mbba_\mbq$ can, in the particular case in section \ref{ss:liftcon}, be ``lifted'' to a categorical connection $\mbba_{\mbp}$ on $\mbp$ such that the pushforward of $\mbba_{\mbp}$ to $\mbq$ is the original connection $\mbba_\mbq$.

We  use  pushforwards to construct an extension of  the notion  of  gauge transformation  of connections.  
 In classical bundle theory, a global gauge transformation is specified by a smooth map $\theta: P\to G$ that is equivariant in the sense that
  \begin{equation}
\theta_{pg}=g^{-1}\theta_pg,
 \end{equation}
 for all $p\in P$ and $g\in G$. 
  In Proposition \ref{P:diffgt} we show that a categorical gauge transformation  on $\mbp^{A, {\rm dec}}$ is specified by a pair $(\theta, \Lambda^H)$, where   $\theta:P\to G$  is as above, and $\Lambda^H$ is an $L(H)$-valued smooth $1$-form  on $P$ that satisfies the equivariance condition
  \begin{equation}
  \Lambda^H_{pg}(vg)=\alpha(g^{-1})\Lambda^H_p(v),
  \end{equation}
  for all $p\in P$, $v\in T_pP$, and $g\in P$. 
  
    A {\em categorical bundle morphism} $\Theta: \mbp^{A, {\rm dec}}\to \mbp^{A, {\rm dec}}$ is a functor that preserves the categorical bundle structure.  In Proposition \ref{P:diffgt} we find the detailed structure of such bundle morphisms.  The morphism $\Theta$   induces a connection on the classical bundle $P$, by keeping track of what happens to the horizontal path $\ovgam$. Thus the traditional connection $A$ gives rise, through this process, to a new connection, and it is this generalized gauge transformation that we introduce and study in section \ref{s:catclassggg}.

   In our concluding result, Proposition \ref{P:gengaugept}, we show that the action of the categorical gauge transformation on categorical connections leads to the following transformation of the classical gauge field $A$:
  \begin{equation}\label{E:introAthtauHL}
  A\mapsto {\rm Ad}(\theta)A -(d\theta)\theta^{-1}+\tau \Lambda^H.
  \end{equation}
In  \cite{CLSgg2018}*{equation (1.2)} we used a different approach, in terms of local trivializations of bundles, to obtain a version of this result, with all forms pulled down to the base manifold. This transformation law (\ref{E:introAthtauHL}) is also superficially similar to the gauge transformation law in higher gauge theories obtained by Wang  \cite{Wang2H}*{equation (1.2)}), within a different framework.

  \subsection{Background in Higher Gauge Theory}\label{ss:ow} Parallel-transport over path spaces have been studied in both the mathematics and physics literature.   We shall mention just a few other works, though there are is now quite a substantial body of literature on different approaches to higher gauge theories. Among the early works that directly or indirectly influenced the study of higher gauge theories is Migdal's work \cite{Mig1980}, wherein a loop-space formulation of quantum chromodynamics was used. Gross \cite{Gross1985} developed a mathematically precise theory of connections over path spaces and derived results for Yang-Mills theory using this framework. 
  Alvarez et al. \cite{Alva1998}, and later \cite{ Alva2009}, studied the problem of finding conserved quantities in integrable field-theoretic systems. Here they considered parallel transport over higher-dimensional geometric objects, and used multiple higher forms, beyond the usual 1-form, for such parallel transport processes. 
  Pfeiffer \cite{Pfei2003} and Girelli et al.  \cite{GirPf2004} used category-theoretic methods.  Other works with a heavier category-theoretic  focus    include Baez et al. \cites{BaezSch,BaezHuerta2011, BaezW}, Martins et al. \cites{MarPick2010, MarPickGray,MarPickSurf }, Parzygnat \cite{Parz2015}, and Sati et al. \cite{SatSchSt}. Higher gauge transformations have been studied within other frameworks by Breen and Messing \cite{BreenMess2005}, Schreiber et al.  \cites{SchreibWalLocal, SchrWalGerb13}, Soncini and Zucchini \cite{SZ}, Waldorf \cites{WalSurfHol10, Waldorf2016, Waldorf2017}, and Wang \cites{Wang3G, Wang3R, Wang2H}. In our approach we don't use \v{C}ech cohomology or gerbes, and  we don't make any use of local trivializations.  
  
  Our approach to categorical principal bundles, developed in \cites{CLS2geom, CLSpdg}, is {\em driven more by geometry than category theory}. This approach, including   gauge transformations as developed in this paper, can be extended to higher path/surface spaces, but in this paper we have focused on path space categories.

 \section{Categorical principal bundles}\label{s:cpb}
 
 In this section we introduce notation, notions, and basic results about categorical groups, path space categories, and categorical principal bundles.
 
 \subsection{A note on notation} We will use a convenient but nonstandard convention of displaying maps or morphisms from right to left rather than left to right. Thus  
  \begin{equation}\label{E:mapconv}
\begin{tikzcd}
q &p \ar[l,bend right=70,swap,"\gamma"] 
\end{tikzcd}
\end{equation}
is a map (or morphism or path)  with domain (source or initial point) $p$ and codomain (target or terminal point) $q$. The advantage of this display convention is that a composition  $\delta\circ\gamma$ is displayed
in the same order as $\delta $ and $\gamma$ appear in $\delta\circ\gamma$ :
  \begin{equation*}
\begin{tikzcd}
r& q \ar[l,bend right=70,swap,"\delta"]  &p \ar[l,bend right=70,swap,"\gamma"] 
\end{tikzcd}
\end{equation*}
This convention has been used extensively by Parzygnat (for example, in \cite{Parz2015}).

 \subsection{Categorical groups} A categorical group  is a category $\mbg$ along with a functor
 $$\mbg\times\mbg\to\mbg$$
 that makes both $\Obj(\mbg)$ and $\Mor(\mbg)$ groups. A categorical Lie group is a categorical group $\mbg$ for which $\Obj(\mbg)$ and $\Mor(\mbg)$ are both Lie groups, the source and target maps 
 $$s, t:\Mor(\mbg)\to\Obj(\mbg)$$
  are smooth and so is the identity-assigning morphism 
 $$\Obj(\mbg)\to\Mor(\mbg):a\mapsto 1_a.$$
 Associated to a categorical group $\mbg$ is a crossed module $(G, H,\alpha,\tau)$, where  
 \begin{equation}
 \begin{split}
 G &=\Obj(\mbg)\\
 H&= \ker s:\Mor(\mbg)\to \Obj(G).
 \end{split}
 \end{equation}
Thus any element of $H$ is a morphism $e\to x$ for some $x\in G$, with $e$ being the identity in $G$. The homomorphism
 $$\tau:H\to G$$
 is just the target map $t$ restricted to $H$, and 
 $$\alpha:G\times H\to H:(g,h)\mapsto \alpha_g(h)$$
 is given by
 $$\alpha_g(h)=1_gh1_{g^{-1}}.$$
 The categorical group $\mbg$ can be reconstructed from $(G, H,\alpha,\tau)$ by taking $\mbg$ to have object group $G$ and morphism group the semi-direct product $H\rtimes_{\alpha}G$.  Henceforth, 
 $$\hbox{\em we will write $(h,g)\in H\rtimes_{\alpha}G$ as $hg$;}$$
  (this  notation carries  a slight risk of confusion but is very convenient). In particular, we identify $g\in G$ with $(e,g)\in H\rtimes_{\alpha}G$ and $h\in H$ with $(h,e)$. Then
 \begin{equation}
 \alpha_g(h)=ghg^{-1}.
 \end{equation}
 We note the Peiffer identities \cite{Peif}:
\begin{equation}\label{E:Peiffer}
\begin{split}
\tau\bigl(\alpha_g(h)\bigr) &= g   \tau(h)  g^{-1}\\
\alpha_{\tau(h)}(h') &=hh'h^{-1}
\end{split}
\end{equation}
for all $g\in G$ and $h\in H$.

With $(h,g)\in H\rtimes_{\alpha}G$ viewed as an element of $\Mor(\mbg)$, the source and targets are
\begin{equation}\label{E:sthg}
s(h,g)=g,\qquad\hbox{and}\qquad t(h,g)=\tau(h)g.
\end{equation}
 Composition of morphisms, viewed as an operation on $H\times G$, is given by
 \begin{equation}
 (h_2, g_2)\circ (h_1, g_1) =(h_2h_1, g_1),
 \end{equation}
 where $g_2=\tau(h_1)g_1$ for the composition to be meaningful.

The categorical group   $\mbg$ is a categorical Lie group if and only if $G$ and $H$ are Lie groups and $\alpha$ and $\tau$ are smooth. 

\subsection{The categorical group $\mbgd$}  For any group $G$, let $\mbgd$ be the categorical group whose objects are the elements of $G$ and for which there is a unique morphism $g_0\to g_1$ for any $g_0, g_1\in G$. Following the notational convention in (\ref{E:mapconv}), we display this unique morphism as
 $$g_1\leftarrow g_0,$$
 where $g_0$ is the source of the morphism and $g_1$ is the target. In the crossed module notation $H\rtimes_{\alpha}G$, the group $H$ is the same as $G$, and the target map is
 $$g_1=t(k, g_0)= kg_0,$$
 for $k\in H=G$. Thus
 \begin{equation}\label{E:GdotdotHG}
 g_1\mapsfrom g_0 \qquad\hbox{corresponds to $(g_1g_0^{-1}, g_0)\in G\rtimes_{\alpha}G$.}
 \end{equation}

 Thus the object group of $\mbgd$ is $G$ (not to be confused with the category whose only object is $G$), and whose morphism group is $G\times G$.
 If $\mbg$ is a categorical group whose object group is $G$ then we have the functor
 \begin{equation}\label{E:defSGdd}
 S:\mbg\to\mbgd,
 \end{equation}
 which is just the identity map $G\to G$ on objects, and takes any $\phi\in\Mor(\mbg)$ to the morphism
\begin{equation}\label{E:defSGddst}
t(\phi)\mapsfrom s(\phi) \end{equation}
 in $\Mor(\mbgd)$. It is readily checked that $S$ is indeed a functor, and, moreover, it is a homorphism of groups at the object and at the morphism levels. If $\mbg$ is a categorical Lie group then so is $\mbgd$ in the obvious way, and $S$ is  smooth both at the object and at the morphism levels.

 \subsection{Smooth  spaces}  We will not need any details concerning smooth structures on path spaces but we note here some minimal background. For our purposes it is convenient to use the framework of {\em diffeological spaces}, introduced by  Souriau \cite{Sour1980} and discussed further by several authors \cites{IZ2013, BH2011, Lau2006}; however, we will use the term ``smooth space'', which is used by Baez and Hoffnung \cite{BH2011} in a broader sense. There are several other approaches to smooth structures, such as the one by Fr\"ohlicher \cite{Froh1982}; Batubenge et al. \cite{batubenge2017diffeological} and Stacey \cite{stacey2010comparative} provide overviews and comparisons of different approaches to smoothness.
 
 Very briefly, we take a smooth space to be a non-empty set $X$ along with a ``diffeology'', a family $S_X$ of maps $U\to X$,  called {\em plots}, with $U$ running over all open subsets of all finite-dimensional spaces $\mbr^n$ with $n\geq 0$, such that: (i) all maps from the one-point space $\mbr^0$ to $X$ are in $S_X$; (ii) if $\phi:U\to X$ is in $S_X$ and if $g:V\to U$ is $C^\infty$, where $V$ is an open subset of some $\mbr^n$, then $\phi\circ g\in S_X$; (iii) if $\phi:U\to X$ is such that the restriction of $\phi$ to every member of an open covering of $U$ is in $S_X$ then $\phi\in S_X$.
 
 If $X$ is a smooth space and $X_0$ is a non-empty subset of $X$, then a diffeology on $X_0$ is obtained by taking as plots all the plots $U\to X$ whose images lie in $X_0$.  In particular, a closed interval $[a,b]$ is a smooth space. 
  
   If $X$ and $Y$ are smooth spaces then a map $F:X\to Y$ is said to be smooth if $F\circ \phi\in S_Y$ for all $\phi\in S_X$.  Thus the condition is that for any plot of $X$ the composition with $F$ is a plot of $Y$.
   
   A surjective map $h:X\to Y$ pushes forward a diffeology $S_X$ on $X$ to a diffeology $h_*S_X$ on $Y$, with $h_*S_X$ consisting of all maps $h\circ\phi$ with $\phi$ running over $S_X$.
 
 \subsection{Path spaces}\label{ss:pathsp}  Let $a<b$ be real numbers, and $C^\infty_0([a,b];X)$ the set of all smooth maps $[a,b]\to X$, where $X$ is a smooth space, that are constant near $a$ and near $b$.   We define a plot for $C^\infty_0([a,b];X)$ to be  a smooth variation of paths on $X$ in the following sense. Let $U$ be any nonempty open subset of some $\mbr^m$, with $m\geq 0$. Consider a map of the form
 $$U\to C^\infty_0([a,b];X): u\mapsto \phi_u,$$
 such that
 $$U\times [a,b]\to X: (u,t)\mapsto \phi_u(t)$$
 is smooth, and there is an $\epsilon>0$ such that, for each $u\in U$, the path $\phi_u$ is constant on $[a, a+\epsilon)$ and on $(b-\epsilon, b]$.   We take all such maps as the plots specifying a diffeology on $C^\infty_0([a,b]; X)$. Varying $a$ and $b$, gives a smooth space that is the disjoint union
 \begin{equation}
 \mcp_1(X)=\cup_{a, b\in\mbr, a<b} C^\infty_0([a,b]; X).
 \end{equation}
 On $\mcp_1(X)$ there is an action of $\mbr$ by time-translation: if $u\in\mbr$  and $\gamma\in  C^\infty_0([a,b]; X)$ then we have the path 
 $$\gamma_{+u}\in  C^\infty_0([a-u,b-u]; X): t\mapsto \gamma(t+u).$$
 Then there is a natural surjection of  $\mcp_1(X)$  onto the quotient space  $\mcp_1(X)/\mbr$. The plots on  $\mcp_1(X)$  composed with this projection give plots on $\mcp_1(X)/\mbr$, and make the latter into  a smooth space.
 
Given a path $\delta:[a,b]\to X$ and a path $\gamma:[b,c]\to X$ we can form a composite path
\begin{equation}
    \gamma\ast\delta:[a,c]\to X: u\mapsto \begin{cases} \delta(u) &\hbox{if $t\in [a,b]$;} \\
    \gamma(u) & \hbox{if $u\in [b,c]$.}
    \end{cases}
\end{equation}

The following result shows that composition of paths is a smooth operation.
 
 \begin{prop}\label{P:composmooth}
 Let $a, b, c\in\mbr$ with $a<b<c$, and let $\mck$ be the subset of $C^\infty_0([b,c]; X)\times C^\infty_0([a,b]; X)$ consisting of all pairs $(\gamma, \delta)$ for which $\gamma(b)=\delta(b)$. Then the composite map
 \begin{equation}
K:  \mck\to C^\infty_0([a,c]; X):(\gamma, \delta)\mapsto \gamma\ast\delta
 \end{equation}
 is smooth.
 \end{prop}
 \begin{proof} A plot   of $\mck$ is a plot of $C^\infty_0([b,c]; X)\times C^\infty_0([a,b]; X)$ that takes values in the subset $\mck$. Thus this plot is of the form
 $$(\phi, \psi)$$
 where there is a non-empty open subset $U$ of $\mbr^m$, for some $m\geq 0$, and $\phi: U\to  C^\infty_0([b,c]; X)$ and $\psi: U\to C^\infty_0([a,b]; X)$  are smooth. Then  $K\circ (\phi, \psi)$ maps $u\in U$ to $\kappa_u=\phi_u\circ\psi_u:[a,c]\to X$. Now $\kappa_u$ is smooth on $[b, c]$, being the same as $\phi_u$ on that interval, and also smooth on $[a, b]$; furthermore, it is constant near $b$. It follows, by using property (iii) of plots, that $\kappa_u$ is smooth on $[a,c]$. Moreover, there is an $\epsilon_1>0$ such that each $\phi_u$ is constant within an $\epsilon_1$-neighborhood of the points $b$ and $c$, and there is an $\epsilon_2>0$ such that each $\psi_u$ is constant within an $\epsilon_2$-neighborhood of  $a$ and $b$. Hence,   $\kappa_u$ is constant within an $\epsilon$-neighborhood of   $a$ and $c$, where $\epsilon=\min\{\epsilon_1, \epsilon_2\}$. Thus $K\circ (\phi, \psi)$  is a plot of $C^\infty_0([a,c]; X)$. This means that $K$ is smooth.
 \end{proof}

  \subsection{Categorical path  spaces}   From a smooth space $X$ we can construct a category ${\mathbb P}_1(X)$ as follows.  The object set of   ${\mathbb P}_1(X)$ is $X$, and the morphism set is $\mcp_1(X)/\mbr$. Let us look at the morphisms in more detail.  The morphisms of ${\mathbb P}_1(X)$ arise from smooth   paths 
 $[a,b]\to X:u\mapsto \gamma(u)$, for any $a, b\in\mbr$ with $a\leq b$,  constant near the endpoints $a$ and $b$, with two such paths identified  if one is obtained from the other by a constant translation of the parameter $u$.  A morphism's source is the initial point, which we often denote $\gamma_0$,  and the target is the terminal point, which we denote $\gamma_1$.  Composition of morphisms corresponds to composition of paths, with the first path terminating at the source of the second.  The set of morphisms $\Mor({\mathbb P}_1(X))$ has a natural smooth space structure, and the source and target maps
 $$s, t:  \Mor({\mathbb P}_1(X)) \to X$$
  are smooth.
  
   It is sometimes convenient to identify paths that are reparametrizations of each other in some way, but we will not use this convention. We also don't use the ``reversal'' of a path as an operation on morphisms.

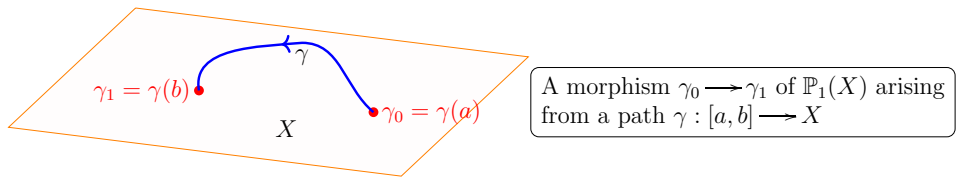
\begin{figure}[ht]

 \tdplotsetmaincoords{70}{110}
 \resizebox{5in}{!}{%
\begin{tikzpicture}[scale=.8,tdplot_main_coords]
   
    \def\x{.5}
   
    \filldraw[
        draw=orange,%
        fill=purple!1,%
    ]          (-1,8,0)
            -- (7,8,0)--(7,17,0)
 -- (-1,17,0)--(-1,8,0);

                \fill [red] (6,12,1) circle[radius=3pt] node[anchor=east] {$\gamma_1=\gamma(b)$};
                
                 \fill [red] (6,16,1) circle[radius=3pt] node[anchor=west] {$\gamma_0=\gamma(a)$};

  \coordinate [label=above: {$\gamma$}] (gam) at (5,14,1.3);

                \coordinate [ label=above: {$X$}] (B) at (6,14,0);

                   \draw[middlearrow={<}, very thick, blue] (6,12,1) to [out=100,in=185]  (5,14,2) to [out=0, in =140] (6,16,1);
                   
                   \node [mybox, right=45,below=12] (box) at (5,22,3){%
    \begin{minipage}{0.55\textwidth}
 A morphism $\gamma_0\to\gamma_1$ of ${\mathbb P}_1(X)$ arising from a path $\gamma:[a,b]\to X$
    \end{minipage}
};

\end{tikzpicture}

}

\caption{The category ${\mathbb P}_1(X)$}
    \label{F:decbun}
    \end{figure}

If $f:X\to Y$ is a smooth map between smooth spaces then $\gamma\mapsto f\circ\gamma$ induces a smooth map 
$${\mathbb P}_1(f):{\mathbb P}_1(X)\to{\mathbb P}_1(Y).$$
In fact, ${\mathbb P}_1$ is a functor from the category of smooth spaces and smooth maps into itself. 
Lastly, let us note that ${\mathbb P}_1$ can be composed with itself multiple times to yield ``higher'' path spaces. There is an alternative, technically easier, way to work with higher path spaces, by viewing them as being obtained from smooth maps $[a_1,b_1]\times\ldots [a_k,b_k]\to X$ that are suitably constant  near the boundary of the domain. We mention this  for cultural context; we will not work with such higher path spaces in this paper.

By a {\em categorical space} we will mean a category for which both object set and morphism set are equipped with smooth space structures such that the following maps are smooth: (i) source and target maps; (ii)  the identity assigning map $a\mapsto 1_a$; (iii) the composition of morphisms, defined on the set of all composable pairs of morphisms.

 There is a special case that we use frequently, for which it is convenient to use a simpler notation. For a smooth manifold $M$ we denote the path space category ${\mathbb P}_1(M)$ by $\mbm$. 
 
 \subsection{Backtrack equivalence}\label{ss:back} There are several reasonable choices for the path space category. One we have used before \cite{CLS2geom}  involves identifying paths that are the same except for some pieces that are backtracked. More precisely, for a path $\gamma:\thinspace  [a,b]\to X$ let ${\gamma_{-1}}:\thinspace  [b,b+b-a]\to X$ be given by ${\gamma_{-1}}(t)=\gamma(2b-t)$. It is reasonable to identify the composite ${\gamma_{-1}}\ast\gamma :\thinspace [a, 2b-a]\to X$ with the constant path at $\gamma(a)$. Next we identify two paths that differ by a finite number of compositions of the type ${\gamma_{-1}}\ast\gamma$. We call this {\em backtrack equivalence}. With this equivalence, the paths with a fixed initial point form a group under composition. In Singer \cite{Sing1995} a principal bundle is constructed informally for which this group of loops based at a fixed point serves as the structure group of a principal bundle.

\subsection{The categorical bundle $\mbpd$} Let
 $$\pi:P\to M$$
 be a principal $G$-bundle, where $G$ is a Lie group. Specifically, $P$ and $M$ are smooth spaces, $\pi$ is a surjective submersion, and there is a smooth free right action of $G$ on $P$:
 $$P\times G\to P:(p,g)\mapsto pg=R_gp$$
 which preserves the fibers of $\pi$.  We will construct categorical spaces from this bundle. Intuitively, the category $\mbpd$ will have the points of $P$ as objects, and morphisms are paths on $M$   connecting the projections on $M$ of the source and the target.


\begin{figure}[ht]

 \tdplotsetmaincoords{70}{110}
\begin{tikzpicture}[scale=.9,tdplot_main_coords]
    \def\x{.5}
   
    \filldraw[
        draw=purple,%
        fill=purple!1,%
    ]          (-1,8,-.5)
            -- (7,8,-.5)--(7,17,-.5)
 -- (-1,17,-.5)--(-1,8,-.5);

         \coordinate [  label=above: {$s({\ovgam})= p_0$}] (sovg) at (6,10.9,3.8);

       
    \fill [ ] (6,12,4) circle[ radius=3pt] node[anchor=south] { 
    };
       \fill (6,16,3.5) circle[radius=2pt] node[anchor=south] {$t({\ovgam})=p_1$};

           \draw[dashed] (6,16, 3.5) --(6,16,0);
           
              \draw[dashed] (6,12,4) --(6,12,0);
              
                \fill [ ] (6,12,0) circle[radius=3pt] node[anchor=west] {$\gamma_0$};
                
                 \fill [ ] (6,16,0) circle[radius=3pt] node[anchor=west] {$\gamma_1$};

  \coordinate [label=above: {$\gamma$}] (gam) at (5,14,-.25);

              \coordinate [label=above: {$ {\ovgam}=(p_1, p_0; \gamma) $}] (ovgam) at (5,13.5, 1.8);

                 \coordinate [ label=above: {$M$}] (M) at (6,17,1);
                
                 \coordinate [  label=above: {$P$}] (P) at (4,17,3);
                 
                   \draw[middlearrow={>}, very thick, red ] (6,12,0) to [out=90,in=185]  (5,12,0) to [out=0, in =180] (6,16,0);

\end{tikzpicture}

\caption{A morphism $\ovgam =(p_1, p_0; \gamma)$ of $\mbpd$.}
    \label{F:Pdotdot}
    \end{figure}

 More precisely, we define $\mbpd$ to be the category whose object set is $P$ and whose morphisms are of the form
 $$(p_1,p_0; \gamma)\in P\times P\times \Mor(\mbm),$$
 with $\gamma$ having source $\pi(p_0)$ and target $\pi(p_1)$:
 
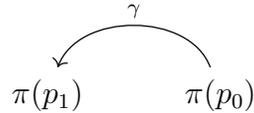
\begin{equation*}
\begin{tikzcd}
\pi(p_1) & \pi(p_0) \ar[l,bend right=70,swap,"\gamma"] 
\end{tikzcd}
\end{equation*}
Source and targets are given by
\begin{equation}\label{E:Pdotsst}
s(p_1, p_0; \gamma)=p_0 \qquad\hbox{and}\qquad t(p_1, p_0; \gamma)=p_1.
\end{equation}
Composition is given by
\begin{equation}\label{E:Pdotscompo}
(p_2, p_1; \delta)\circ (p_1, p_0; \gamma)= (p_2, p_0; \delta\circ\gamma).
\end{equation}
The identity morphism at $p$ is $(p, p; 1_{\pi(p)}),$ where $1_{u}$ is the point-path at $u$.

   The categorical group $\mbgd $ has a  categorical right action on $\mbpd$ given on objects by the action of $G$ on $P$ and on morphisms by
 \begin{equation}
 (p_1,p_0;\gamma)(g_1\mapsfrom g_0)=(p_1g_1, p_0g_0;\gamma).
 \end{equation}
 We have the projection functor
 $$\pi:\mbpd\to\mbm,$$
 given on objects by the bundle projection $\pi:P\to M$ and on morphisms by
 $$\pi(p_1, p_0;\gamma)=\gamma.$$
 
 \subsection{The decorated  categorical bundle  $\mbp^{A, {\rm dec}}$}\label{ss:catdec}
 Now consider a connection $A$ on a principal $G$-bundle $\pi:P\to M$, and let $\mbg$ be a categorical Lie group with associated Lie crossed module $(G, H,\alpha,\tau)$.  From this we can construct a categorical principal $\mbg$-bundle  that we call {\em decorated} bundle and denote 
\begin{equation}\label{E:decbun}
\pi:\mbp^{A, {\rm dec}}\to\mbm.
\end{equation}
   The object space of $\mbp^{A, {\rm dec}}$ is $P$. A morphism of  $\mbp^{A, {\rm dec}}$ is to be thought of as an $A$-horizontal path on $P$ equipped with a decorating element drawn from $H$. More precisely, a morphism of  $\mbp^{A, {\rm dec}}$ is of the form
   $$\ovgam =(\tlg, h),$$
   where $\tlg$ is a morphism of $\mbp$ coming from an  $A$-horizontal path on $P$, smooth and constant near its initial and terminal points, and $h\in H$.  Source and target maps are defined by:
 \begin{equation}\label{E:stMorPA}
 \begin{split}
 s(\ovgam; h) &=s(\tlg)\\
  t(\ovgam; h) &=t(\tlg)\tau(h),
 \end{split}
 \end{equation}
 where $s(\tlg)=\tlg_0$ is the initial point of $\tlg$ and $t(\tlg)=\tlg_1$ is the terminal point. 
 We call  (\ref{E:decbun}) the {\em decorated bundle} corresponding to the bundle $\pi:P\to M$ and connection $A$.  
 

\begin{figure}[ht]

 \tdplotsetmaincoords{70}{110}
\begin{tikzpicture}[scale=.9,tdplot_main_coords]
    \def\x{.5}
   
    \filldraw[
        draw=purple,%
        fill=purple!1,%
    ]          (-1,8,-.1)
            -- (7,8,-.1)--(7,17,-.1)
 -- (-1,17,-.1)--(-1,8,-.1);
       
       \draw[middlearrow={>}, thick,  brown, dashed ] (6,12,4) to (6,16,5);
       
         \coordinate [  label=above: {$s({\ovgam},h)={\ovgam}_0$}] (sovg) at (6,10.5,3.8);

         \coordinate [  label=above: {$({\ovgam}, h)$}] (sovgh) at (6,14,4.5);

       \draw[middlearrow={>}, very thick, brown ] (6,12,4) to [out=90,in=185]  (5,12,2) to [out=0, in =180] (6,16,3.5);
       
    \fill [ ] (6,12,4) circle[ radius=3pt] node[anchor=south] { 
    };
        \fill (6,16,3.5) circle[radius=2pt] node[anchor=south] {${\ovgam}_1$};
        
           \fill [ ] (6,16,5) circle[radius=3pt] node[anchor=west] {$t({\ovgam},h)={\ovgam}_1\tau(h)$};
           
           \draw[dashed] (6,16,5) --(6,16,0);
           
              \draw[dashed] (6,12,4) --(6,12,0);
              
                \fill [ ] (6,12,0) circle[radius=3pt] node[anchor=west] {$\gamma_0$};
                
                 \fill [ ] (6,16,0) circle[radius=3pt] node[anchor=west] {$\gamma_1$};

  \coordinate [label=above: {$\gamma$}] (gam) at (5,14,-.25);

              \coordinate [label=above: {$ {\ovgam} $}] (ovgam) at (5,14,2.5);

                 \coordinate [ label=above: {$M$}] (M) at (6,17,1);
                
                 \coordinate [  label=above: {$P$}] (P) at (4,17,3);
                 
                   \draw[middlearrow={>}, very thick, red ] (6,12,0) to [out=90,in=185]  (5,12,0) to [out=0, in =180] (6,16,0);

\end{tikzpicture}

\caption{The decorated bundle, showing the source and target of a morphism $(\ovgam, h)$. Only  $\bigl(\ovgam_1, \ovgam_0; \gamma\bigr)$ is needed about $\ovgam$.}
    \label{F:decbun2}
    \end{figure}
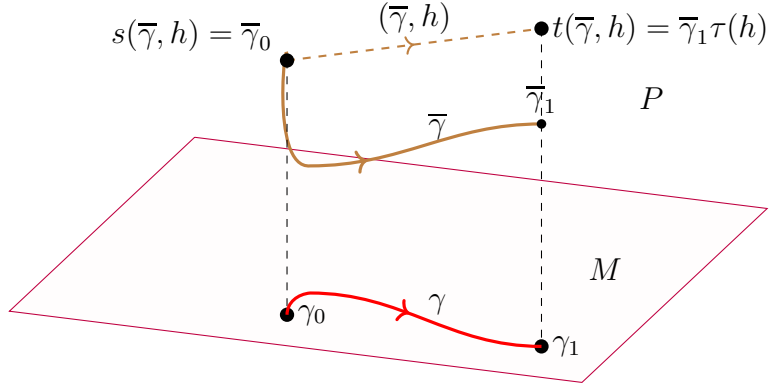


 The categorical group $\mbg$ acts on $\mbp^{A, {\rm dec}}$ on objects by the action of $G$ on $P$ and on morphisms by:
 \begin{equation}\label{E:Gactsdec}
 (\tlg; h)h'g' =({\tlg} g'; \gamma;  {g'}^{-1}hh'g').
 \end{equation}
 Here it is useful to recall that the notation $hg$ is really a short form of $(h,g)$. 
 For $ (\tld; h_2),  (\tlg; h_1) \in \Mor(\mbp^{A, {\rm dec}})$, with 
 \begin{equation}
 \tld_0 ={\tlg}_1\tau(h_1),
 \end{equation}
 the  composition of morphisms is defined by
  \begin{equation}\label{E:compo}
 (\tld; h_2)\circ (\tlg; h_1)= \bigl({\tld} \tau(h_1)^{-1}\circ \tlg; h_1h_2\bigr).
 \end{equation}
 The source of either side is $\tlg_0$ and the target of either side is $\tld_1\tau(h_2)$. Moreover, the identity morphism at $p$ is $(1_p; e)$, where $1_p$ is the constant point-path at $p$.

 We   check the behavior of compositions under the action of the categorical group. Consider the composition (\ref{E:compo}) under the action of a composition $(h'_2g'_2)\circ(h'_1g'_1)$; for the latter to be meaningful we have $g'_2=\tau(h'_1)g'_1$. With this, we have:
 
 \begin{equation}\label{E:functcompo}
 \begin{split}
 &  (\tld; h_2)(h'_2g'_2)\circ (\tlg; h_1)(h'_1g'_1)\\
 &=\bigl({\tld} g'_2;  {g'}_2^{-1}h_2h'_2g'_2\bigr)\circ\bigl({\tlg} g'_1; {g'}_1^{-1}h_1h'_1g'_1\bigr)\\
 &=\bigl({\tld}\tau(h_1)^{-1}g'_1\circ {\tlg}g'_1
 ;  \, {g'}_1^{-1}h_1h_2h'_2h'_1g'_1\bigr),
 \end{split}
 \end{equation}
 which agrees with
  \begin{equation}\label{E:functcompo2}
 \begin{split}
 &  \bigl( (\tld; h_2)\circ (\tlg; h_1)\bigr)\bigl((h'_2g'_2)\circ (h'_1g'_1)\bigr) \\
 &=\bigl({\tld}\tau(h_1)^{-1}\circ \tlg; \, h_1h_2\bigr)(h'_2h'_1g'_1)\\
 &= \Bigl(\bigl({\tld}\tau(h_1)^{-1} \circ {\tlg}\bigr) g'_1; \, {g'}_1^{-1}h_1h_2h'_2h'_1g'_1 \Bigr)\\
 &= \bigl({\tld}\tau(h_1)^{-1}g'_1\circ {\tlg}g'_1; \, {g'}_1^{-1}h_1h_2h'_2h'_1g'_1 \bigr).
  \end{split}
 \end{equation} 
 This means that compositions of morphisms  in $\mbg$ and $\mbp^{A, {\rm dec}}$ commute with the action of $\Mor(\mbg)$ on  $\Mor(\mbp^{A, {\rm dec}})$.

  
  \subsection{Categorical principal bundles} We have discussed $\mbpd$ and $\mbp^{A, {\rm dec}}$. These are both examples of the general notion of a categorical principal bundle as defined in our earlier work  \cite{CLS2geom}.  Let $\mbg$ be a categorical Lie group, and let $\mbm$ be a categorical space; we have in mind the usual case where $\mbm$ arises from a manifold $M$. A {\em categorical principal bundle} with {\em structure categorical   group} $\mbg$ is comprised of   categorical spaces  $\mbp$ and $\mbm$, a functor 
   $$\pi: \mbp \to{\mbm} $$
   that is smooth and surjective  both on the level of objects and on the level of morphisms, along with a functor
      $$\mbp\times\mbg\to\mbp$$
      that is a free smooth  right action   both on objects and on morphisms, such that $\pi(pg)=\pi(p)$ for all objects/morphisms $p$ of $\mbp$ and all objects/morphisms $g$ of $\mbg$.   In practice we are only concerned with the case where    $\mbg$ is a categorical Lie group, $\Obj(\mbp)$ and $\Obj({\mbm})$ are smooth manifolds,  and the object bundle
      $$\Obj(\mbp)\to\Obj({\mbm})$$
      is a principal $G$-bundle, where $G=\Obj(\mbg)$.

 \subsection{The   functor $\mbbs$}\label{ss:functmbs}  Let $\mbp\to\mbm$ be a categorical principal $\mbg$-bundle, and $\mbpd\to\mbm$ the categorical $\mbgd$-bundle discussed earlier, obtained from the object principal $G$-bundle $\pi:\thinspace P\to M$. Let
  \begin{equation}\label{E:defSfunct1}
{\mbbs}:\thinspace  \mbp \to \mbpd,
 \end{equation}
 be given on objects by $p\mapsto p$ and on morphisms by
 \begin{equation}\label{E:defSfunct12}
{\mbbs}({\tilde\gamma})= ({\tilde\gamma}_1, {\tilde\gamma}_0; \pi({\tilde\gamma})),
\end{equation}
where the subscripts $0$ and $1$ signify source and target, respectively. It is readily verified that this is a functor (commutes with source and targets, respects compositions, and maps identities to identities). Moreover, for any $\phi\in\Mor(\mbg)$, we also have
\begin{equation}\label{E:mbsfunct}
\mbs({\tilde\gamma}\phi)= \bigl({\tilde\gamma}_1\phi_1, {\tilde\gamma}_0\phi_0; \pi({\tilde\gamma})\bigr)= \mbs({\tilde\gamma})S(\phi),
\end{equation}
where the first equality holds because of the functorial nature of the action of $\mbg$ on $\mbp$ and the second equality is verified from the definition of $S$ given in  (\ref{E:defSGdd}) and (\ref{E:defSGddst}).

\subsection{The   functor $\mbbs$ for decorated bundles}   Now we specialize to the case where $\mbp\to \mbm$ is the usual decorated bundle $ \mbp^{A, {\rm dec}}\to \mbm$, with $A$ being a connection on the underlying object bundle $P\to M$, as discussed in section \ref{ss:catdec}. We recall that the  $ \mbp^{A, {\rm dec}}$  and  $\mbpd$ both have $P$ as object space, but morphisms of $\mbpd$ are of the form $(p_1, p_0;\gamma)$, where $\gamma\in\Mor(\mbm)$ runs from source $\pi(p_0)$ to target $\pi(p_1)$, while morphisms of $ \mbp^{A, {\rm dec}}$ are of the form $(\tlg; h)$, where now $\tlg$ is an $A$-horizontal morphism of $\mbp$ and $h$ is a general element of $H$.  Then we have the functor
  \begin{equation}\label{E:defSfunct}
{\mbbs} : \mbp^{A, {\rm dec}}\to \mbpd,
 \end{equation}
 given on objects by $p\mapsto p$ and on morphisms by
 \begin{equation}\label{E:defSfunct2}
{\mbbs}(\tlg; h)= (\tlg_1\tau(h), \tlg_0;\gamma),
 \end{equation}
 where $\tlg_0=s(\tlg)$ and $\tlg_1=t(\tlg)$.
 Let us verify  for this case the properties of $\mbs$ noted in the general context earlier. 
 
 \begin{prop}\label{P:Pifunct} The assignment ${\mbbs} $ given above is a functor. Moreover, it is a  morphism of categorical principal bundles in the following sense:
 \begin{equation}\label{E:Sfunct}
 \begin{split}
 {\mbbs}(pg) &={\mbbs} (p)g\qquad\hbox{for all $(p,g)\in P\times G$,}\\
 {\mbbs}\bigl((\tlg; h)h'g'\bigr) &={\mbbs}(\tlg; h)S(h', g') \end{split}
 \end{equation}
 for all $(\tlg; h)\in \Mor\left(\mbp^{A, {\rm dec}}\right)$ and $(h',g')\in H\rtimes_{\alpha}G$, and $S$ is as given in (\ref{E:defSGdd}) and (\ref{E:defSGddst}).
  \end{prop}
 \begin{proof} From (\ref{E:defSfunct2}) it is readily seen that ${\mbbs} $ maps sources and targets properly. Moreover, for compositions we have:
 \begin{equation}
 \begin{split}
  {\mbbs}\bigl( (\tld; h_2)\circ (\tlg; h_1)\bigr) 
 &={\mbbs}({\tld}\tau(h_1)^{-1}\circ\tlg; h_1h_2)\\
 &=\bigl({\tld}_1\tau(h_2),  \tlg_0;\, \delta \circ\gamma),
 \end{split}
 \end{equation}
 which agrees with
  \begin{equation}
 \begin{split}
  {\mbbs}(\tld; h_2)\circ {\mbbs}(\tlg; h_1)  
 &= \bigl(\tld_1\tau(h_2), \tld_0; \delta\bigr)\circ\bigl(\tlg_1\tau(h_1), \tlg_0 ;\gamma\bigr)\\
 &=\bigl(\tld_1\tau(h_2), \tlg_0; \delta\circ\gamma),
 \end{split}
 \end{equation}
 where we used the composition law (\ref{E:Pdotscompo}).

 The first equation in (\ref{E:Sfunct}) is immediate from the definition of ${\mbbs} $ in (\ref{E:defSfunct}) acting on objects. Next, using the definition (\ref{E:defSfunct2}) of ${\mbbs} $ on morphisms and the action given by (\ref{E:Gactsdec}), we have
 \begin{equation}
 \begin{split}
  {\mbbs}\bigl((\tlg; h)h'g'\bigr) &={\mbbs}({\tlg}g';  {g'}^{-1}hh'g')\\
  &=\bigl({\tlg}_1\tau(hh')g', \tlg_0g' ; \gamma\bigr),
 \end{split} 
 \end{equation}
 and
 \begin{equation}
 \begin{split}
 {\mbbs}(\tlg; h){S}(h'g') &=
  \bigl(\tlg_1\tau(h), \tlg_0; \gamma) \bigl(\tau(h')g'\mapsfrom g'\bigr)\\
  &=\bigl(\tlg_1\tau(hh')g', \tlg_0g';\gamma\bigr).
 \end{split}
 \end{equation}
 Thus ${\mbbs}$ satisfies (\ref{E:Sfunct}). 
  \end{proof}

 \section{Connections on Categorical Bundles}\label{s:ConCat}
 
 Let $\mbg$ be a categorical Lie group with associated Lie crossed module $(G, H, \alpha, \tau)$. We will now look at a counterpart of some of the previously discussed constructions, but with the traditional principal bundle $\pi:P\to M$, with classical connection $A$, replaced by a categorical bundle $\pi:\mbp\to\mbm$ with a categorical analog of a classical connection.
 
\subsection{Categorical connections}\label{ss:catcon}  A connection $A$ on a principal $G$-bundle specifies a special path $\tlg^A_p$ on $P$, starting at any given point $p$ on the fiber over $\gamma_0$, and is called the $A$-horizontal lift of $\gamma$ starting at $p$.  This generalizes then readily to categorical bundles. A {\em categorical connection}  $\mbba$ on a categorical principal $\mbg$-bundle $\pi:\mbp\to\mbm$ assigns to each $\gamma\in\Mor(\mbm)$ and each $p\in \Obj(\mbp)$ with $\pi(p)=s(\gamma)$, a morphism, called the {\em horizontal lift}, 
 $$\tau_{\mbba}(\gamma; p)\in\Mor(\mbp)$$ with source $p$ and whose $\pi$-projection is $\gamma$, such that the following conditions hold:
 \begin{itemize}
 \item[(CC1)] If  $\gamma=1_u$, the identity   at $u=\pi(p)\in \Obj(\mbm)$, then  $\tau_{\mbba}(\gamma; p)=1_p$;
 \item[(CC2)] $\tau_{\mbba}(\gamma; pg)= \tau_{\mbba}(\gamma; p)1_g$ for all $g\in G$;
 \item[(CC3)] If $\gamma, \delta \in\Mor(\mbm)$ are such that the composite $\delta\circ\gamma$ is defined then  the horizontal lift $\tau_{\mbba}(\delta\circ\gamma; p)$ is the composite of the horizontal lift $\tau_{\mbba}(\gamma ;p)$ followed by the horizontal lift of $\delta$:
\begin{equation}\label{E:taucompo}
\tau_{\mbba}(\delta\circ\gamma; p)= \tau_{\mbba}\bigl(\delta; t\bigl(\tau_{\mbba}(\gamma ;p)\bigr)\bigr)\circ \tau_{\mbba}(\gamma;p).
\end{equation}
 \end{itemize}
 We also require that $(\gamma; p)\to\tau_{\mbba}(\gamma; p)$ be smooth, where, of course, $p\in \Obj(\mbp)$ and $\gamma\in\Mor(\mbm)$  are such that $\pi(p)=s(\gamma)$,
 
 \subsection{The standard example} Let $A$ be a connection on a classical principal $G$-bundle $\pi:P\to M$, where $G=\Obj(\mbg)$ is the object group of a categorical Lie group $\mbg$ associated to the crossed module $(G, H,\alpha,\tau)$. Then we can construct a categorical connection $\mbbad$ on the categorical $\mbgd$-bundle $\pi:\mbpd\to\mbm$ by setting
 \begin{equation}
 \tau_{\mbbad}(\gamma; p)= (q, p; \gamma),
 \end{equation}
 where $q$ is the point  obtained by   parallel transporting $p$ along $\gamma$ by $A$.
 
 Intuitively, every categorical connection on $\mbpd$ arises in this way from a classical connection on $P$.  
 
 \subsection{The horizontal bundle   $\mbp^{\mbba }$}  With setting as above, by a {\em horizontal morphism}  we shall mean a morphism of the form   $\tau_{\mbba}(\gamma; p)$. Property (CC3) implies that the composition of horizontal lifts is horizontal. Thus we have a category $\mbp^{\mbba}$, whose object set is $\Obj(\mbp)$ and whose morphisms are  all the horizontal morphisms. The categorical group involved for this bundle has object group $G$ and the only morphisms are the identity morphisms $1_g$ for all $g\in G$.
  
\subsection{The decorated bundle   $\mbp^{\mbba, {\rm dec}}$}\label{ss:catdec2} In section \ref{ss:catdec} we saw how a classical connection $A$ on a principal $G$-bundle $\pi:P\to M$, along with a categorical Lie group $\mbg$ whose object group is $G$, lead to a categorical $\mbg$-bundle $\pi:\mbp^{A, {\rm dec}}\to \mbm$. Here we shall see how this process generalizes to a categorical connection on a categorical principal bundle.   Let $\mbg$ and $\mbg_1$ be categorical Lie groups, with  
$$\Obj(\mbg_1)=\Obj(\mbg)=G.$$
We have in mind the case  where $\mbg_1$ is $\mbgd$.  As before we  have taken $(G, H, \alpha, \tau)$ to be the Lie crossed module associated to $\mbg$.

Let  $\mbba$  be a categorical connection on a  categorical principal $\mbg_1$-bundle $\pi:\mbp \to\mbm$. Thus to each $\gamma\in\Mor(\mbm)$ and $p\in \Obj(\mbp)$ there is the associated  horizontal lift
 $ \tau_{\mbba}(\gamma; p)$. The {\em decorated} categorical  principal $\mbg$-bundle 
 \begin{equation}
 \mbp^{\mbba, {\rm dec}},
 \end{equation}
has
 \begin{equation}
 \Obj(\mbp^{\mbba, {\rm dec}})=P, \qquad \Mor(\mbp^{\mbba, {\rm dec}})= \Mor(\mbp^{\mbba})\times H.
 \end{equation}
 Source and targets are given by
 \begin{equation}
 s(\ovgam;h)=\ovgam_0, \qquad t(\ovgam; h)=\ovgam_1 \tau(h),
 \end{equation}
 where, as always,   the subscripts $0$ and $1$ signify source and target, respectively.
 We have the functor 
  \begin{equation}
 \mbbs: \mbp^{\mbba, {\rm dec}}\to \mbpd,
 \end{equation}
 given on objects by $p\mapsto p$ and on morphisms by
 \begin{equation}
 \mbbs(\ovgam; h)= \bigl( \ovgam_1\tau(h), \ovgam_0; \pi(\ovgam)\bigr).
 \end{equation}

 This is the more abstract form of the functor $\mbbs$ introduced in (\ref{E:Sfunct}) in the context of the
decorated categorical bundle arising from a classical principal $G$-bundle with connection for $A$.
 
 \subsection{Lifting a connection}\label{ss:liftcon} Our goal here and in the next subsection is to transfer connections from one bundle to another. First, here, we see how
 a connection on $\mbpd$ can be lifted to a connection on $ \mbp^{\mbba, {\rm dec}}$, by decorating each horizontal morphism in $\mbpd$ with the identity $e\in H$.

\begin{prop}\label{P:liftcon}
With notation as above, let $\mbba$ be a categorical connection on  the categorical $\mbgd$-bundle $ \mbpd $. For any $\gamma \in \Mor(\mbm) $ and $p\in P$ on the fiber over $s(\gamma)$, let
\begin{equation}\label{E:deftau1}
\tau_{\mbba_d}(\gamma;p)= \bigl(\tau_{\mbba}( \gamma; p), e\bigr) \in \Mor( \mbp^{\mbba, {\rm dec}}).
\end{equation}
Then $\mbba_d$ is a categorical connection on  $ \mbp^{\mbba, {\rm dec}}$.
\end{prop}
\begin{proof}  First we note that  $\bigl(\tau_{\mbba}(\gamma; p), e\bigr)$ is indeed in $ \Mor( \mbp^{\mbba, {\rm dec}})$.  The condition (CC1) is readily verified. Condition (CC2) follows by applying the definition (\ref{E:Gactsdec}) of the action of $H\rtimes_{\alpha}G$ on $ \Mor( \mbp^{\mbba, {\rm dec}})$. Lastly, (CC3) follows by using the composition specification (\ref{E:compo}). Smoothness of $(\gamma;p)\mapsto  \bigl(\tau_{\mbba}( \gamma; p), e\bigr)$ follows from smoothness of $\tau_{\mbba}$.
\end{proof}

\subsection{Pushing forward a connection}\label{ss:projcon}  Let $\pi:\mbp\to\mbm$ be a categorical $\mbg$-bundle and $\pi:\mbpd\to\mbm$ the corresponding categorical $\mbgd$-bundle. We have then the functor
\begin{equation}\label{E:defSfunct1b}
\mbs:\mbp\to\mbpd,
\end{equation}
introduced in (\ref{E:defSfunct1}).
Now suppose $\mbba$ is a categorical connection on $\mbp\to\mbm$. We construct a categorical connection $\mbbad$ on $\mbpd\to\mbm$ as follows. For $\gamma\in\Mor(\mbm)$ and $p\in \Obj(\mbp)$ on the fiber above the source $\gamma_0$, we set
\begin{equation}\label{E:deftauadot}
\tau_{\mbbad}(\gamma; p)=\mbs\bigl(\tau_{\mbba}(\gamma; p)\bigr).
\end{equation}
Because of the properties of $\mbs$ this assignment defines a categorical connection $\mbbad$.

Specializing this pushforward process  to a connection on $ \mbp^{\mbba, {\rm dec}}$ produces a connection on  $\mbpd$, which we verify in the following result.

\begin{prop}\label{P:connproj} With notation as above, suppose $\mbba_1$ is a categorical connection on $ \mbp^{\mbba, {\rm dec}}$.  If $\tau_{\mbba_1}(\gamma;p)=\bigl(\tlg, h\bigr)$, where $\gamma\in\Mor(\mbp)$ and $p\in \Obj(\mbp)$, with $\pi(p)=s(\gamma)$, we set  
\begin{equation}\label{E:deftauA1dots}
\tau_{\mbbad_1}(\gamma;p)=\bigl(q\tau(h), p; \gamma\bigr)\in \Mor(\mbpd),
\end{equation}
where   $q=t(\tlg)$.
 Then $\mbbad_1$ is a categorical connection on $\mbpd$.
\end{prop}
\begin{proof} For notational convenience let us write $\mbba_2$ for $\mbbad_1$.

 For (CC1) we note that if $\gamma=1_u$, where $u=\pi(p)$, then 
\begin{equation}\label{E:tau11p}
\tau_{\mbba_1}(1_{u}; p)= (1_p, e),
\end{equation}
 and so
\begin{equation}\label{E:tau21p}
\tau_{\mbba_2}(1_u; p)=\bigl(p, p; 1_u\bigr).
\end{equation}
  
 Next, for (CC2),  if $\tau_{\mbba_1}(\gamma;p)=(\tlg, h)$, where $\tlg$ runs from $p$ to $q$,  then
 \begin{equation}
 (\tlg, h)1_g =({\tlg} 1_g, g^{-1}hg)
 \end{equation}
 and so, noting that
 ${\tlg}1_g$ runs from $pg$ to $pq$, we have
 \begin{equation}
\tau_{\mbba_2}(\gamma; pg) = (qg\,g^{-1}\tau(h)g, pg; \gamma)= (q\tau(h)g, pg; \gamma).
\end{equation}
This agrees with:
\begin{equation}
\tau_{\mbba_2}(\gamma; p)1_g = \bigl(q\tau(h), p; \gamma\bigr)1_g= \bigl(q\tau(h)g, pg; \gamma\bigr).
\end{equation}
Thus (CC2) holds.

Now, for (CC3), we consider a composite $\delta\circ\gamma$ and $p\in P$ with $\pi(p)=u=s(\gamma)$. Let
\begin{equation}
\tau_{\mbba_1}(\gamma; p)=(\tlg, h_1)\qquad\hbox{and}\qquad \tau_{\mbba_1}(\delta; q)=(\tld,  h_2),
\end{equation}
where $q=t(\tlg)$.  

Then
\begin{equation}\label{E:tauA2prcomp}
\begin{split}
\tau_{\mbba_1}(\delta\circ\gamma; p) &=\tau_{\mbba_1}\bigl(\delta; t\bigl(\tau_{\mbba_1}(\gamma; p)\bigr)\bigr)\circ\tau_{\mbba_1}(\gamma,p)\\  
&= \tau_{\mbba_1}\bigl(\delta; q\tau(h_1)\bigr)\circ (\tlg, h_1)\\
&=\bigl(  {\tld} \tau(h_1), \tau(h_1)^{-1}h_2\tau(h_1)\bigr) \circ (\tlg, h_1)\\
&= \bigl(  {\tld} \tau(h_1),  h_1^{-1}h_2h_1\bigr) \circ (\tlg, h_1)\\
&=\bigl(\tld\circ\tlg, h_2h_1)
\end{split}
\end{equation}
where we used the second of the Peiffer identities (\ref{E:Peiffer}). Hence,
\begin{equation}\label{E:tauA2dgp}
\tau_{\mbba_2}(\delta\circ\gamma; p) = \bigl(r\tau(h_2h_1), p; \delta\circ\gamma).
\end{equation}
On the other hand,

\begin{equation}
\begin{split}
&\tau_{\mbba_2}\bigl({\delta}; t\bigl(\tau_{\mbba_2}(\gamma;p)\bigr)\bigr)\circ\tau_{\mbba_2}(\gamma,p)\\
&=\tau_{\mbba_2}\bigl({\delta};    q\tau(h_1)    \bigr)\circ \bigl(p\tau(h_1), p; \gamma\bigr)
\\
 &=\bigl(r\tau(h_1) \tau(h_1)^{-1}\tau(h_2)\tau(h_1), q\tau(h_1);{\delta}\bigr)\circ \bigl(p\tau(h_1), p; \gamma\bigr)\\
 &=\bigl(r\tau(h_2h_1), q\tau(h_1);{\delta}\bigr)\circ\bigl(p_1\tau(h_1), p; \gamma\bigr)\\
&=\bigl(r\tau(h_2h_1), p; {\delta}\circ\gamma\bigr),
\end{split}
\end{equation}
which agrees   with  the right side in (\ref{E:tauA2prcomp}). Thus (CC3) holds. 

Smoothness of $\tau_{\mbba_2}$ follows from smoothness of $\tau_{\mbba_1}$ and of $\tau$.
\end{proof}

\section{Pushforwards of Categorical Connections}\label{s:pfcon}

In this section we present a more abstract construction of the pushforward discussed in section \ref{ss:projcon}.

Let us briefly recall how classical connections can be pushed forward  from one bundle to another. For a detailed and more general account we refer to Kobayashi  and Nomizu 
\cite{KobNomI}*{section II.6}; there is also a pullback process for classical bundles that we will not discuss here. In the categorical bundles framework we also have a pullback in the special situation discussed in section \ref{ss:liftcon}.

\subsection{Pushforwards of classical connections} Let 
$G$ and $K$ be Lie groups, and $s:G\to K$ a smooth homomorphism. Now consider a commutative diagram
\begin{equation}\label{E:PQmap}
\begin{tikzcd}[column sep={3em,between origins}]
P\arrow{rr}{S}\arrow{dr}[swap]{\pi_{P}} &&
Q   \arrow{dl}{\pi_{Q}}   \\ 
& M  
\end{tikzcd}
\end{equation}
where $\pi_P:P\to M$ is a principal $G$-bundle and $\pi_Q:Q\to M$ is a principal $K$-bundle, and suppose that the smooth map $S$ is equivariant in the following sense
\begin{equation}\label{E:PQmap2}
S(pg)= S(p)s(g),
\end{equation}
for all $p\in P$ and $g\in G$.   

With this setting there is a way to push forward a connection $A$ on $P$ to a connection $f_*A$ on $Q$ as follows.  Let $q\in Q$ and $v\in T_qQ$. We pick any point $p\in P$ on the fiber over $\pi_Q(q)$; then $q=S(p)k$, for some $k\in K$. Then the horizontal space $\ker (f_*A)_q$ for $f_*A$ at $q$ is $f_*(\ker A_p)k$. More geometrically, {\em $f_*A$-horizontal paths in $Q$ are $K$-translates of the images under $f$ of $A$-horizontal paths in $P$. }


\subsection{Pushforwards of categorical connections}  We now  construct a categorical counterpart of this pushforward process. Let $ \mathbf P$ and $ \mathbf Q$  be, respectively, a $ \mathcal G$- and an $ \mathcal K$- categorical principal  bundle  over the same base category $ \mathbf M$. Here $ \mathcal G$ and a $ \mathcal K$ are categorical Lie groups. Suppose  $ {{\mbs}}: \mathbf P \to \mathbf Q$ and $ {S}: \mathcal G \to \mathcal K$ are functors satisfying
\begin{equation}\label{E:defmbspg1}
\begin{split}
\mbs(pg) &=\mbs(p)S(g)\\
 {{\mbs}}({\tilde \gamma} \phi)&={{\mbs}}({\tilde \gamma} )S(\phi)\\
 \end{split}
 \end{equation}
  for all $p\in\Obj(\mbp)$, $g\in\Obj(\mcg)$, $ \tlg\in \Mor({ \mathbf P})$ and  $ \phi \in \Mor(\mathcal G)$.   We assume also that
  \begin{equation}\label{E:defmbspg2}
 \pi_{\mbq}\bigl(\mbs(p)\bigr)=\pi_{\mbp}(p)  \end{equation}
 for all $p\in\Obj(\mbp)$.  We say that the pair $\mbs$ and $S$, satisfying (\ref{E:defmbspg1}) and (\ref{E:defmbspg2}), is a morphism from the bundle $\mbp$ to the bundle $\mbq$.
  
  A categorical connection is known when all the  horiziontal morphisms are known. Starting from a categorical connection $\tau_{\mbp}$ on $\mbp$ we {\em specify a categorical connection $\mbq$ by requiring that  the $\mbs$-images of  $\tau_{\mbp}$-horizontal morphisms be $\tau_{\mbq}$-horizontal.} The following result ensures that $\tau_{\mbq}$ is a categorical connection.
  
  To avoid notational clutter we will write ${\ovgam}1_{b}$ as ${\ovgam}b$, for $\ovgam$ a morphism in $\mbp$ or $\mbq$ and $b$ an object of $\mcg$ or $\mck$. This is consistent with standard practice for the case where $\ovgam$ corresponds to a horizontal path.

\begin{prop}\label{P:pushfor}  
We use the notation as above. Let $\tau_{\mbp}$ be a categorical connection on $\mbp$, and, for any $\gamma\in \Mor(\mbm)$ and any $q\in \Obj(\mbq)$ on the  $\mbq$-fiber over the source $s(\gamma)$, let
 \begin{equation}\label{E:tauQ}
\tau_{\mathbf Q}(\gamma, q):={{\mbs}}(\tau_{\mathbf P}(\gamma, p)) {k_{p,q}},
\end{equation}
where $p\in\Obj(\mbp)$ is any point on the $\mbp$-fiber over $s(\gamma)$ and $k_{pq}\in \Obj({\mathcal K})$ is specified by requiring that $q=\mbs(p)k_{p,q}$. Then $\tau_{\mbq}$ satisfies conditions (CC1)-(CC3) for categorical connections,
\end{prop}

\begin{proof} First, let us verify that (\ref{E:tauQ}) is independent of the choice of the point $p$. Suppose $a\in \Obj(\mcg)$ is such that $\mbs(pa)=\mbs(p)=q$. Then 
$$S(a)=e,$$
 the identity in the group $\Obj(\mck)$.  We observe that
$$q=\mbs(p)k_{p,q}= \mbs(pa) k_{p,q},$$
and so
\begin{equation}
k_{pa, q}=k_{p,q}.
\end{equation}
Then
\begin{equation}
\begin{split}
{{\mbs}}(\tau_{\mathbf P}(\gamma, pa)) {k_{pa, q}} &={{\mbs}}(\tau_{\mathbf P}(\gamma, p)1_a) {k_{pa, q}} \\
&= {{\mbs}}(\tau_{\mathbf P}(\gamma, p)){S}(1_a) {k_{p, q}}\\
&= {{\mbs}}(\tau_{\mathbf P}(\gamma, p)){}1_{S(a)} {k_{p, q}}\\
&={{\mbs}}(\tau_{\mathbf P}(\gamma, p)){}  {k_{p, q}},
\end{split}
\end{equation}
since $S(a)=e$. Thus (\ref{E:tauQ}) is independent of the choice of $p$.

Next we verify the conditions (CC1-3)  from section \ref{ss:catcon}.  Condition (CC1) is readily verified and we omit the argument.

For (CC2)  consider $\gamma\in\Mor(\mbm)$, $q\in \Obj(\mbq)$ a point lying on the fiber above the source $s(\gamma)$.  We choose any $p\in\Obj(\mbp)$ on the fiber over $s(\gamma)$, and let $k_{p,q}\in \Obj(\mck)$ be such that $q=\mbs(p)k_{p,q}$. Then, for $k\in  \Obj(\mck)$, we have
$$qk=\mbs(p)k_{p,q}k,$$
and so
\begin{equation}
k_{p, qk}=k_{p,q}k.
\end{equation}
Hence
\begin{equation}
\begin{split}
 \tau_{\mathbf Q}(\gamma, qk) &={{\mbs}}(\tau_{\mathbf P}(\gamma, p)) k_{p,q}k\\
& =\tau_{\mathbf Q}(\gamma, q)k.
\end{split}
\end{equation}
Thus, the translate of any $\tau_{\mbq}$-horizontal morphisms by any $k\in \Obj(\mck)$ is $\tau_{\mbq}$-horizontal.

Next, we verify that the condition (CC3)  holds; this  condition means   that compositions of $\tau_{\mbq}$-horizontal morphisms are $\tau_{\mbq}$-horizontal.  By definition of $\tau_{\mbq}$, as given in (\ref{E:tauQ}), a $\tau_{\mbq}$-horizontal morphism is of the form
$$\tau_{\mbq}(\gamma, q)=\mbs\bigl(\tau_{\mbp}(\gamma, p)\bigr)k_{p,q}.$$
Let $\gamma_1, \gamma_2\in\Mor(\mbm)$ and $p_1, p_2\in \Obj(\mbp)$ be such that the composition
\begin{equation}\label{E:Staucompo}
\mbs\bigl(\tau_{\mbp}(\gamma_2, p_2)\bigr)k_{p_2, q_2}\circ  \mbs\bigl(\tau_{\mbp}(\gamma_1, p_1)\bigr)k_{p_1, q_1}
\end{equation}
is defined.  Let $p'_1$ be the terminal point (target) of  $\tau_{\mbp}(\gamma_1, p_1)$.  Since the composition (\ref{E:Staucompo}) is defined we have
\begin{equation}\label{E:Sp2kpq}
\mbs(p_2)k_{p_2,q_2}=\mbs(p'_1)k_{p_1,q_1}.
\end{equation}
Applying $\pi_{\mbq}$ shows that
$$\pi_{\mbp}(p_2)=\pi_{\mbp}(p'_1),$$
and so the composition $\gamma_2\circ\gamma_1$ is defined. Moreover, 
$$p_2=p'_1a,$$
for some $a\in \Obj(\mcg)$. Then, applying $\mbs$ and using   (\ref{E:Sp2kpq}), we see that
\begin{equation}
S(a)k_{p_2,q_2}=k_{p_1,q_1}.
\end{equation}
Using this relation in (\ref{E:Staucompo}) turns it into
\begin{equation}\label{E:Staucompo2}
\bigl[\mbs\bigl(\tau_{\mbp}(\gamma_2, p_2)\bigr) \circ  \mbs\bigl(\tau_{\mbp}(\gamma_1, p_1)\bigr)S(a) \bigr]k_{p_2, q_2}.
\end{equation}
The term within $[\cdots]$ is equal to
\begin{equation}
\mbs\bigl(\tau_{\mbp}(\gamma_2, p_2)\bigr) \circ  \mbs\bigl(\tau_{\mbp}(\gamma_1, p)\bigr),
\end{equation}
where
$$p=p_1a.$$
Since $\mbs$ is a functor, we have
\begin{equation}
\begin{split}
\mbs\bigl(\tau_{\mbp}(\gamma_2, p_2)\bigr) \circ  \mbs\bigl(\tau_{\mbp}(\gamma_1, p)\bigr)&=
\mbs\bigl(\tau_{\mbp}(\gamma_2, p_2)\circ \tau_{\mbp}(\gamma_1, p)\bigr)\\
&=\mbs\bigl(\tau_{\mbp}(\gamma_2\circ\gamma_1, p)\bigr).
\end{split}
\end{equation}
Combining this with the expression in (\ref{E:Staucompo2}) shows that the original composition of $\tau_{\mbq}$-horizontal morphisms given in (\ref{E:Staucompo}) is equal to
\begin{equation}
\mbs\bigl(\tau_{\mbp}(\gamma_2\circ\gamma_1, p)\bigr)k_{p_2,q_2}.
\end{equation}
By definition of $\tau_{\mbq}$, as given in (\ref{E:tauQ}), we have
\begin{equation}
\mbs\bigl(\tau_{\mbp}(\gamma_2\circ\gamma_1, p)\bigr)k_{p_2,q_2}=\tau_{\mbq}(\gamma_2\circ\gamma_1, q'_2),
\end{equation}
where $q'_2=\mbs(p)k_{p_2,q_2}$. Thus, the composition of $\tau_\mbq$-horizontal morphisms given in  (\ref{E:Staucompo}) is the $\tau_{\mbq}$-horizontal morphism $\tau_{\mbq}(\gamma_2\circ\gamma_1, q'_2)$. \end{proof}


\section{Gauge Transformations}\label{ss:gt}

 By a {\em gauge transformation} on the categorical $\mbg$-bundle $\pi: \mbp\to\mbm$ we mean a functor
 $$\Theta:\mbp\to\mbp$$
 that commutes with the action of $\mbg$,  satisfies $\pi\circ\Theta=\pi$, and is smooth on both objects and morphisms, $\Theta^{-1}$ exists and is also smooth on objects and morphisms. In \cite{CLSgg2018} we have studied gauge transformations mainly in the context of decorated categorical bundles and twisted product bundles.
 
\subsection{Gauge transformations for objects and morphisms} For any $p\in \Obj(\mbp)$ the point (object) $\Theta(p)$ lies on the same fiber as $p$ and so there is a $\theta_p\in G$ such that
 $$\Theta_p=p\theta_p.$$
 Next, for the same reason, for any $\ovgam\in\Mor(\mbp)$ we have
 $$\Theta(\ovgam)={\ovgam}\theta(\ovgam),$$
 for some $\theta(\ovgam)\in \Mor(\mbg)$. Now since $\Mor(\mbg)\simeq H\rtimes_{\alpha}G$ we can write any element of $\Mor(\mbg)$ as a product of an element of $H$ and an element of $G$, in either order, with both of these groups being viewed as subgroups of  $\Mor(\mbg)$. Moreover, $gh$, with $g\in G$ and $h\in H$, has source $g$ and target $g\tau(h)g^{-1}$.
 
 Thus we can write 
 $$\theta(\ovgam)=h_{\ovgam}g_{\ovgam}$$
 for some $g_{\ovgam}\in G$, which is the source of $\theta(\ovgam)\in\Mor(\mbg)$, and $h_{\ovgam}\in H$. Thus,
 \begin{equation}
 \begin{split}
 \Theta(p) &=p{\theta_p}\\
 \Theta(\ovgam) &={\ovgam}h_{\ovgam}g_{\ovgam}
 \end{split}
 \end{equation}
 for some $g_\ovgam\in G$, $h_{\ovgam}\in H$. The source of $\Theta(\ovgam)$ is then $\ovgam_0g_{\ovgam}$. For this to match $\Theta(\ovgam_0)$, which is $\ovgam_0\theta_{\ovgam_0}$, we have
 \begin{equation}
 g_{\ovgam}=\theta_{\ovgam_0}.
 \end{equation}
 
\subsection{Gauge transformations for morphism compositions} Next consider a composition $\ovdel\circ\ovgam$ of morphisms $\ovdel$ and $\ovgam$ in $\mbp$.  
 
  \begin{equation*}
\begin{tikzcd}
                r & \ar[l,bend right=70,swap,"\ovdel"]         q &p \ar[l,bend right=70,swap,"\ovgam"] 
\end{tikzcd}
\end{equation*}

For $\Theta $ to be a functor, $\Theta $ applied to $\ovdel\circ\ovgam$ must agree with $\Theta (\ovdel)\circ\Theta (\ovgam)$. For the latter we compute
 \begin{equation}
 \begin{split}
 \Theta (\ovdel)\circ\Theta (\ovgam) &= {\ovdel}h_{\ovdel}\theta_q\circ {\ovgam}h_{\ovgam}{\theta_p}\\
 &=(\ovdel\circ\ovgam)\bigl((h_{\ovdel},\theta_q)\circ (h_{\ovgam}, {\theta_p})\bigr)\\
 &=(\ovdel\circ\ovgam)(h_{\ovdel}h_{\ovgam}, {\theta_p}).
 \end{split}
 \end{equation}
 For this to agree with $\Theta (\ovdel\circ\ovgam)$, whose value is given by
 $$\Theta (\ovdel\circ\ovgam)=(\ovdel\circ\ovgam)(h_{\ovdel\circ\ovgam}, {\theta_p}),$$
 the condition is
 \begin{equation}
 h_{\ovdel\circ\ovgam}=h_{\ovdel}h_{\ovgam}.
 \end{equation}
 
 We can summarize this discussion in the following conclusion.
 
 \begin{prop}\label{P:Thetacompo} Let $\Theta :\mbp\to\mbp$ be a gauge transformation for the categorical $\mbg$-bundle $\pi:\mbp\to\mbm$, where the categorical Lie group $\mbg$ has associate Lie crossed module $(G, H, \alpha, \tau)$. Then
 \begin{equation}\label{E:Thetagtpgam}
  \begin{split}
 \Theta(p) &=p{\theta_p}\\
 \Theta(\ovgam) &={\ovgam}h_{\ovgam}g_{\ovgam},
  \end{split}
  \end{equation}
  for all $p\in \Obj(\mbp)$ and $\ovgam\in \Mor(\mbp)$, where ${\theta_p}, g_{\ovgam}\in G$ and $h_{\ovgam}\in H$, and, furthermore,
 \begin{equation}\label{E:gpovgamcompo}
 \begin{split}
 g_{\ovgam} &={\theta_p}\\
  h_{\ovdel\circ\ovgam}&=h_{\ovdel}h_{\ovgam},
  \end{split}
  \end{equation}
  whenever $\ovgam, \ovdel\in \Mor(\mbp)$ are composable morphisms with $p=s(\ovgam)$,
  and
  \begin{equation}\label{E:ghequiv}
  \begin{split}
  \theta_{pa} &=a^{-1}{\theta_p}a \\
  h_{{\ovgam}a} &= a^{-1}h_{\ovgam}a\\
   h_{{\ovgam}b} &= h_{\ovgam}\theta_pb\theta_p^{-1}
  \end{split}
  \end{equation}
  for all $a\in G$ and $b\in H$, with $p$ being $s(\ovgam)$.
  
  Conversely, if $\Theta:\mbp\to\mbp$ is such that   conditions  (\ref{E:Thetagtpgam}), (\ref{E:gpovgamcompo}) and (\ref{E:ghequiv}) hold then $\Theta$ is a functor that commutes with the projection $\pi_{\mbp}$ and with the action of the categorical group $\mbg$.
  \end{prop}
  \begin{proof} We have already established all claims except for (\ref{E:ghequiv}). Next, we have, for $p\in P$ and $a\in G$,
  \begin{equation}
  \Theta(pa)=\Theta(p)a=  paa^{-1}\theta_pa,
  \end{equation}
  which proves the first equation in (\ref{E:ghequiv}).  Using this, we see that for $\ovgam\in\Mor(\mbp)$ with source $p$, and for any $a\in G$, we have
    \begin{equation}
  \Theta({\ovgam}a)=\Theta({\ovgam})a=  {\ovgam}h_{\ovgam}{\theta_p}a={\ovgam}a a^{-1}h_{\ovgam}a{\theta_{pa}},
  \end{equation}
  which implies the second relation in (\ref{E:ghequiv}). Finally, for $b\in H$, noting that
  $$s({\ovgam}b)=s(\ovgam)s(b)=pe=p,$$
   we have
  \begin{equation}
  \Theta({\ovgam}b) = {\ovgam}h_{\ovgam}\theta_{p}b={\ovgam}(h_{\ovgam}\theta_{p}b{\theta_p}^{-1}){\theta_p}.
  \end{equation}
 This  shows that $h_{{\ovgam}b}=h_{\ovgam}{\theta_p}b{\theta_p}^{-1}$, which is the third equation in (\ref{E:ghequiv}). 
    \end{proof}

  \subsection{Gauge transformations at the differential level}\label{ss:gtdiff}  We work now with a classical connection $A$ on a principal $G$-bundle $\pi:P\to M$, and the  decorated categorical principal bundle $\mbp^{A, {\rm dec}}$. Let $(G, H,\alpha, \tau)$ be a Lie crossed module with associate categorical group denoted $\mbg$.  Let $\theta:P\to G$ be a smooth mapping that satisfies
\begin{equation}\label{E:defthetapa}
\theta_{pa}=a^{-1}\theta_pa \qquad\hbox{for all $a\in G$ and $p\in P$.}
\end{equation}
  Let $\Lambda^H$ be a smooth $1$-form on $P$ with values in the Lie algebra $L(H)$ satisfying
  \begin{equation}\label{E:lamdequiv}
  \Lambda^H_{pa}(va)=\alpha(a^{-1})\Lambda^H_p(v),
  \end{equation}
  for all $p\in P$, $a\in G$, and $v\in T_pP$, and vanishes on vertical vectors (that is, vectors in the kernel of $d\pi$).
  
   For any $A$-horizontal path $\tlg:[a,b]\to P$ we  {\em define  
 \begin{equation}\label{E:defhtlg}
 h_{\tlg}=h(b), 
 \end{equation}
   where $h:[a,b]\to H$ is the solution to the differential equation}
  \begin{equation}\label{E:hudiffeq}
 \frac{dh(u)}{du} h(u)^{-1} =-\Lambda^H\bigl(\tlg'(u)\bigr) 
  \end{equation}
  In the following result we show how $\Lambda^H$ and $\theta$ give rise to a gauge transformation $\Theta$ on the categorical bundle $\mbp^{A, {\rm dec}}\to\mbm$.    We denote a typical $A$-horizontal path on $P$ by $\tlg$ or $\tld$.

   \begin{prop}\label{P:diffgt}
   With $\theta$ and $\Lambda^H$, and other notation  as above, let
   \begin{equation}\label{E:defThthLamd}
   \Theta: \mbp^{A, {\rm dec}}\to \mbp^{A, {\rm dec}}
   \end{equation}
   be given  on objects and on morphisms by
    \begin{equation}\label{E:Thetagtpgamh}
  \begin{split}
 \Theta(p) &=p{\theta_p}\\
 \Theta(\ovgam) &={\ovgam}h_{\tlg} \theta_p h,
  \end{split}
  \end{equation}
  for all $p\in \Obj(\mbp)$ and $\ovgam\in \Mor(\mbp^{A, {\rm dec}})$ of the form 
  $$\ovgam =(\tlg, h),$$
  where  $p =s (\tlg)$. Then $\Theta$ is a gauge transformation.
      \end{prop}
      The second condition in (\ref{E:Thetagtpgam}) is motivated by the third equation in (\ref{E:ghequiv}) holds.  It is easier to understand if we write the condition as
      \begin{equation}\label{E:Thetovgamhg}
      \Theta(\ovgam)= {\ovgam}h_{\ovgam}\theta_p,
      \end{equation}
      with
      \begin{equation}\label{E:defhovgam}
      h_{\ovgam}=h_{\tlg}\theta_p h\theta_p^{-1}.
      \end{equation}
      In the conclusion what is missing from saying that $\Theta$ is a gauge transformation is the smoothness.

  \begin{proof} We will verify the conditions given in Proposition \ref{P:Thetacompo} hold.
  
  Consider smooth  $A$-horizontal paths  $\tlg:[a,b]\to P$ and  $\tld:[b,c]\to P$, constant near the initial and terminal points, that are composable; this means, 
  $$\tld(b)=\tlg(b).$$
   On the interval $[a,c]$ let us consider the map
\begin{equation}
f:[a,c]\to H: u\mapsto \begin{cases} h_1(u) \quad\hbox{if $u\in [a,b]$;}\\
  h_2(u)  h_1(b)  \quad\hbox{if $u\in [b,c]$,}
  \end{cases}
  \end{equation}
  where $h_1(\cdot)$ is a solution to (\ref{E:hudiffeq})  for $u\in [a,b]$ with $h_1(a)=e$, and $h_2(\cdot)$ is a solution to (\ref{E:hudiffeq})  for $u\in [b, c]$ with $h_2(b)=e$. Then $f$ solves the equation (\ref{E:hudiffeq}) for $h$, with $f(a)=e$.  Thus
  \begin{equation}
  \frac{df(u)}{du}f(u)^{-1} = -\Lambda^H\bigl((\tld\circ\tlg)'(u)\bigr)
  \end{equation}
  for all $u\in [a,c]$ and $f(u)=e$. Hence
  $$f(c)=h_{\tld\circ\tlg}, $$
 and so
 \begin{equation}
 h_{\tld\circ\tlg}=h_2(c)h_1(b)= h_{\tld} h_{\tlg}.
 \end{equation}
 
 The condition (\ref{E:lamdequiv}) ensures that 
  \begin{equation}\label{E:hovgama}
  h_{{\tlg}a}=\alpha(a^{-1})h_{\tlg} = a^{-1}h_{\tlg}a,
  \end{equation}

   Recall from       (\ref{E:Thetovgamhg}) and  (\ref{E:defhovgam})  that for a morphism  $\ovgam$ of $\mbp^{A, {\rm dec}}$ given by $(\tlg, h)$, with source $p$, we  have defined
   \begin{equation}\label{E:defhovgam2}
   h_{\ovgam}=h_{\tlg}\theta_p h\theta_p^{-1}.
   \end{equation}
For any $b\in H$  we have
   $${\ovgam} b=(\tlg, hb),$$
   and so
    \begin{equation}\label{E:hovgamb}
   h_{{\ovgam}b} =h_{\tlg}\theta_p hb\theta_p^{-1}=h_{\tlg}\theta_p h \theta_p^{-1}\theta_p b\theta_p^{-1}=h_{\ovgam} \theta_p b\theta_p^{-1},
   \end{equation}
   which shows that  the third equation in (\ref{E:ghequiv}) holds. 
   
   Finally, for $a\in G$, we have
   \begin{equation}
   {\ovgam}a=({\tlg} a, a^{-1}ha),
   \end{equation}
   and so 
   \begin{equation}
   \begin{split}
   h_{{\ovgam}a} &= h_{{\tlg}a} \theta_{pa} (a^{-1}ha) \theta_{pa}^{-1} \hskip 1in \hbox{(from the definition (\ref{E:defhovgam}))}
   \\
   &=h_{{\tlg}a} a^{-1}\theta_p h\theta_p^{-1}a\\
   &= a^{-1}h_{\tlg}a   a^{-1}\theta_p h\theta_p^{-1}a \hskip 1in \hbox{(using (\ref{E:hovgama}))}\\
   &= a^{-1}h_{\tlg} \theta_ph\theta_p^{-1}a\\
   &=a^{-1}h_{\ovgam}a  \hskip 1.8in \hbox{(using (\ref{E:defhovgam2}))}.
   \end{split}
   \end{equation}
   Thus we have verified that the functions $p\mapsto \theta_p$ and $\ovgam\mapsto h_{\ovgam}$ satisfy all the properties,  other than smoothness, listed in Proposition \ref{P:Thetacompo} associated with a categorical gauge transformation.  Smoothness of $\Theta$   on objects follows from the smoothness of the mapping $\theta:P\to G$. Smoothness of $\Theta$ on morphisms follows from the fact that if $(u, v)\mapsto \tlg(u; v)\in P$ is smooth  for $u\in [a,b]$ and $v$ running over some $\mbr^k$, then the solution $u\mapsto h(u; v)$  to the differential equation (\ref{E:hudiffeq}), with $\tlg'(u; v)$ instead of $\tlg'(u)$ on the right side,  depends smoothly also on the parameter $v$.  \end{proof}
 
 \section{Categorical gauge transformations of classical connections}\label{s:catclassggg}
 
 Let $A$ be a connection on a principal $G$-bundle $\pi:P\to M$. We then have the categorical connection $A^{\bullet\bullet}$ on the categorical $\mbgd$-bundle $\pi:\mbpd\to\mbm$ and we have the decorated categorical $\mbg$-bundle $\pi:\mbp^{A, {\rm dec}}\to \mbm$. Let $\gamma\in\Mor(\mbm)$. Then the $A$-horizontal lift  $\tlg^A_p$ of $\gamma$, with initial point $p$, gives a morphism 
 $${\ovgam}_p^A=(\tlg^A_p; e)\in \Mor\bigl(\mbp^{A, {\rm dec}}\bigr).$$
 Applying a gauge transformation
  \begin{equation}
  \Theta:\mbp^{A, {\rm dec}}\to \mbp^{A, {\rm dec}}
  \end{equation}
   to the morphism  ${\ovgam}_p^A$ yields the morphism
 \begin{equation}
\Theta({\ovgam}_p^A) =\bigl(\tlg\theta_p; h_{\tilde\gamma}\bigr),
\end{equation}
where, for notational simplicity, we are writing $\tlg^A_p$ simply as $\tlg$, and we are using other notation as before.
This then pushes down to the morphism
\begin{equation}
\bigl(qg_q\tau(h_{\tilde\gamma}), p{\theta_p}; \gamma\bigr)\in \Mor(\mbpd),
\end{equation}
where $q$ is the terminal point of $\tlg^A_p$. 
Let us assume that the new categorical connection on $\mbpd$ also arises from a classical connection on $P$.
Thus the resulting  connection  on $P$ parallel transports the point $p$ along $\gamma$ to the point $q\theta_q\tau(h_{\tilde\gamma}){\theta_p}^{-1}$.

Thus the new horizontal path is
\begin{equation}\label{E:newhor1}
t\mapsto \ovgam(t)= \tlg_t \theta_{\tlg_t}\tau(h_{\tlg_t})\theta_p^{-1},
\end{equation}
where, of course, $\tlg$ is $A$-horizontal, with initial point $p$. 
 
\subsection{Parallel transport with respect to a shifted connection}  Suppose $A$  is a connection form on $\pi:P\to M$, and an $L(G)$-valued $1$-form $C$ on $P$ that satisfies the ${\rm Ad}_G$-equivariance
$$C_{pg}(v_pg)={\rm Ad}(g)^{-1}C_p(v_p)\qquad\hbox{for all $p\in P$,}$$
and vanishes  on vertical vectors:
$$C_p(V_p)=0\qquad\hbox{for all $V_p\in\ker d\pi_p$ and all $p\in P$.}$$
Then $A+C$ is also a connection form. We note how  $(A+C)$-horizontal paths are given by suitable right-translations of  $A$-horizontal paths.

\begin{lem}\label{L: partransshift}
Let $A$ be a connection on a principal $G$-bundle $\pi:P\to M$, and $C$ an ${\rm Ad}_G$-equivariant $1$-form on $P$ with values in $L(G)$ that vanishes on vertical vectors. Suppose $[a,b]\to P: t\mapsto {\tlg}(t)$ is an $A$-horizontal path. Then $[a,b]\to P:t\mapsto {\tlg}(t)\xi(t)$ is horizontal with respect to the connection $A+C$, if $t\mapsto \xi(t)$ satisfies the differential equation
\begin{equation}\label{E:difeqxi}
\xi'(t)\xi(t)^{-1}= -C\bigl(\tlg'(t)\bigr)\qquad \hbox{for all $t\in [a,b]$.}
\end{equation}
\end{lem}
\begin{proof}  We apply $A+C$ to $t\mapsto {\tlg(t)}\xi(t)$ to compute:
\begin{equation}
\begin{split}
(A+C)\bigl(\tlg'(t)\xi(t)+ \tlg(t)\xi'(t)\bigr) &= 0 +\xi(t)^{-1}\xi'(t) +\xi(t)^{-1}C\bigl(\tlg'(t)\bigr)\xi(t)+0\\
&=\xi(t)^{-1}\left[\xi'(t)\xi(t)^{-1}+ C\bigl(\tlg'(t)\bigr) \right]\xi(t).
\end{split}
\end{equation}
Thus the path $[a,b]\to P:t\mapsto \tlg(t)\xi(t)$ is $(A+C)$-horizontal if and only if equation (\ref{E:difeqxi}) holds. \end{proof}

\subsection{Traditional gauge transformation of $A$} Let  $\phi:P\to P$ be a traditional gauge transformation; this means that it is a  smooth $G$-equivariant map for which $\pi\circ \phi=\pi$.  Since $p$ and $\phi(p)$ are on the same $\pi$-fiber, there is a unique $\theta(p) \in G$ such that
\begin{equation}\label{E:phip}
\phi(p)=p\theta(p)\qquad\hbox{for all $p\in P$.}
\end{equation}
Local triviality can be used to show that $\theta$ is smooth, and $G$-equivariance of $\phi$ is equivalent to the condition
\begin{equation}
\theta(pg)=g^{-1}\theta(p)g\qquad\hbox{for all $p\in P$ and $g\in G$.}
\end{equation}
(Proofs may be found in any standard text on bundle theory such as \cite{KobNomI}.) Then
\begin{equation}
(\phi^{-1})^*A = {\rm Ad}(\theta_p^{-1})A_p  -(d\theta|_p)\theta_p^{-1}.
\end{equation}
If $\tlg$ is $A$-horizontal then $\phi\circ \tlg$ is $(\phi^{-1})^*A$-horizontal:
\begin{equation}\label{E:phiinvA}
(\phi^{-1})^*A\bigl((\phi\circ\tlg)'(t)\bigr)= A\bigl(\phi^{-1}_*(\phi\circ\tlg)'(t)\bigr)=A\bigl(\tlg'(t)\bigr)=0.
\end{equation}

\subsection{Horizontal paths and generalized gauge transformations} In view of (\ref{E:phiinvA}) and Lemma \ref{L: partransshift}, if $\tlg$ is an $A$-horizontal path on $P$, and $\theta:P\to G$ is as above, associated to a gauge transformation $\phi:P\to P$, then the path
\begin{equation}
[a,b]\to P: t\mapsto \tlg(t)\theta\bigl(\tlg(t)\bigr)\xi(t)
\end{equation}
is horizontal with respect to the connection
\begin{equation}\label{E:AthtauHL}
{\rm Ad}(\theta)A -(d\theta)\theta^{-1} +\tau\Lambda^H,
\end{equation}
where
\begin{equation}
\xi'(t)\xi(t)^{-1} =- \tau\Lambda^H\bigl(\tlg'(t)\bigr)
\end{equation}
for all $t\in [a,b]$. We can write $\xi$ as
\begin{equation}
\xi(t)= \tau\bigl(h_{\tlg}(t)\bigr),
\end{equation}
where $t\mapsto h_{\tlg}(t)$ solves
\begin{equation}
\frac{dh_{\tlg}(t)}{dt}h_{\tlg}(t)^{-1} =-\Lambda^H\bigl(\tlg'(t)\bigr).
\end{equation}
Thus, we have established the following result.

\begin{prop}\label{P:gengaugept}
With notation as above, and $\tlg:[a,b]\to P$ the $A$-horizontal path on $P$ with initial point $p$, the path
\begin{equation}
[a,b]\to P: t\mapsto  \tlg(t)\theta_{\tlg(t)} \tau\bigl(h_{\tlg}(t)\bigr) \theta(p)^{-1}
\end{equation}
has initial point $p$ and
is horizontal with respect to the connection 
\begin{equation}\label{E:ggA}
{\rm Ad}(\theta)A -(d\theta)\theta^{-1} +\tau\Lambda^H.
\end{equation}
\end{prop}
Looking back at (\ref{E:newhor1}) we conclude that the connection form induced by  the gauge transformation $\Theta$ on $\mbp^{A, {\rm dec}} $ is given by (\ref{E:ggA}).

 \subsection{Concluding remarks}
 
 In this work we have constructed pushforwards of categorical connections. Applying this to a decorated categorical bundle, whose ingredients are a classical principal $G$-bundle $\pi:P\to M$ with a connection  and a second Lie structure group $H$, intertwined with $G$ in a specific way, we obtain a transformation of $A$ of the form $\phi^*A+\tau\Lambda^H$, where $\phi:P\to P$ is a classical gauge transformation and $\tau\Lambda^H$ arises from a particular type of $L(H)$-valued $1$ form on $P$. 
 
 We have applied the pushforward method only to bundles over path spaces. However, the same ideas should be applicable to higher path spaces, and the corresponding gauge transformations would involve not only the classical transformation $\phi$ and the   $1$-form $\Lambda^H$ but also higher-order forms with values in other Lie algebras.
           
        {\bf{ Acknowledgments.} }   Chatterjee acknowledges  research support from SERB, DST, Government of India grant MTR/2018/000528.   Lahiri thanks the S.N. Bose National
Centre for a travel grant to attend a research workshop. Sengupta thanks Arthur Parzygnat for discussions. The authors thank the University of Connecticut for research support for a workshop related to this work.

\bibliography{CLS}{}

\begin{bibdiv}
\begin{biblist}

\bib{Alva1998}{article}{
      author={Alvarez, Orlando},
      author={Ferreira, Luiz~A.},
      author={Guill\'{e}n, J.~S\'{a}nchez},
       title={A new approach to integrable theories in any dimension},
        date={1998},
        ISSN={0550-3213},
     journal={Nuclear Phys. B},
      volume={529},
      number={3},
       pages={689\ndash 736},
}

\bib{Alva2009}{article}{
      author={Alvarez, Orlando},
      author={Ferreira, Luiz~A.},
      author={S\'{a}nchez-Guill\'{e}n, J.},
       title={Integrable theories and loop spaces: fundamentals, applications
  and new developments},
        date={2009},
        ISSN={0217-751X},
     journal={Internat. J. Modern Phys. A},
      volume={24},
      number={10},
       pages={1825\ndash 1888},
}

\bib{BH2011}{article}{
      author={Baez, John~C.},
      author={Hoffnung, Alexander~E.},
       title={Convenient categories of smooth spaces},
        date={2011},
        ISSN={0002-9947},
     journal={Trans. Amer. Math. Soc.},
      volume={363},
      number={11},
       pages={5789\ndash 5825},
}

\bib{BaezHuerta2011}{article}{
      author={Baez, John~C.},
      author={Huerta, John},
       title={An invitation to higher gauge theory},
        date={2011},
        ISSN={0001-7701},
     journal={Gen. Relativity Gravitation},
      volume={43},
      number={9},
       pages={2335\ndash 2392},
}

\bib{BaezSch}{incollection}{
      author={Baez, John~C.},
      author={Schreiber, Urs},
       title={Higher gauge theory},
        date={2007},
   booktitle={Categories in algebra, geometry and mathematical physics},
      series={Contemp. Math.},
      volume={431},
   publisher={Amer. Math. Soc., Providence, RI},
       pages={7\ndash 30},
}

\bib{BaezW}{article}{
      author={Baez, John~C.},
      author={Wise, Derek~K.},
       title={Teleparallel gravity as a higher gauge theory},
        date={2015},
        ISSN={0010-3616},
     journal={Comm. Math. Phys.},
      volume={333},
      number={1},
       pages={153\ndash 186},
}

\bib{batubenge2017diffeological}{misc}{
      author={Batubenge, Augustin},
      author={Iglesias-Zemmour, Patrick},
      author={Karshon, Yael},
      author={Watts, Jordan},
       title={Diffeological, {F}r\"{o}licher, and {D}ifferential {S}paces},
        date={2017},
}

\bib{BreenMess2005}{article}{
      author={Breen, Lawrence},
      author={Messing, William},
       title={Differential geometry of gerbes},
        date={2005},
        ISSN={0001-8708},
     journal={Adv. Math.},
      volume={198},
      number={2},
       pages={732\ndash 846},
         url={http://dx.doi.org/10.1016/j.aim.2005.06.014},
      review={\MR{2183393}},
}

\bib{CLS2geom}{article}{
      author={Chatterjee, Saikat},
      author={Lahiri, Amitabha},
      author={Sengupta, Ambar~N.},
       title={Path space connections and categorical geometry},
        date={2014},
        ISSN={0393-0440},
     journal={J. Geom. Phys.},
      volume={75},
       pages={129\ndash 161},
}

\bib{CLSpdg}{article}{
      author={Chatterjee, Saikat},
      author={Lahiri, Amitabha},
      author={Sengupta, Ambar~N.},
       title={Connections on decorated path space bundles},
        date={2017},
        ISSN={0393-0440},
     journal={J. Geom. Phys.},
      volume={112},
       pages={147\ndash 174},
}

\bib{CLSgg2018}{article}{
      author={Chatterjee, Saikat},
      author={Lahiri, Amitabha},
      author={Sengupta, Ambar~N.},
       title={Gauge transformations for categorical bundles},
        date={2018},
        ISSN={0393-0440},
     journal={J. Geom. Phys.},
      volume={133},
       pages={219\ndash 241},
}

\bib{Froh1982}{incollection}{
      author={Fr\"{o}licher, Alfred},
       title={Smooth structures},
        date={1982},
   booktitle={Category theory ({G}ummersbach, 1981)},
      series={Lecture Notes in Math.},
      volume={962},
   publisher={Springer, Berlin-New York},
       pages={69\ndash 81},
      review={\MR{682945}},
}

\bib{GirPf2004}{article}{
      author={Girelli, Florian},
      author={Pfeiffer, Hendryk},
       title={Higher gauge theory---differential versus integral formulation},
        date={2004},
        ISSN={0022-2488},
     journal={J. Math. Phys.},
      volume={45},
      number={10},
       pages={3949\ndash 3971},
}

\bib{Gross1985}{article}{
      author={Gross, Leonard},
       title={A {P}oincar\'{e} lemma for connection forms},
        date={1985},
        ISSN={0022-1236},
     journal={J. Funct. Anal.},
      volume={63},
      number={1},
       pages={1\ndash 46},
         url={https://doi-org.libezp.lib.lsu.edu/10.1016/0022-1236(85)90096-5},
      review={\MR{795515}},
}

\bib{IZ2013}{book}{
      author={Iglesias-Zemmour, Patrick},
       title={Diffeology},
      series={Mathematical Surveys and Monographs},
   publisher={American Mathematical Society, Providence, RI},
        date={2013},
      volume={185},
        ISBN={978-0-8218-9131-5},
}

\bib{KobNomI}{book}{
      author={Kobayashi, Shoshichi},
      author={Nomizu, Katsumi},
       title={Foundations of differential geometry. {V}ol. {I}},
      series={Wiley Classics Library},
   publisher={John Wiley \& Sons, Inc., New York},
        date={1996},
        ISBN={0-471-15733-3},
        note={Reprint of the 1963 original, A Wiley-Interscience Publication},
}

\bib{Lau2006}{article}{
      author={Laubinger, Martin},
       title={Diffeological spaces},
        date={2006},
        ISSN={0716-0917},
     journal={Proyecciones},
      volume={25},
      number={2},
       pages={151\ndash 178},
}

\bib{MarPick2010}{article}{
      author={Martins, Jo\~ao~Faria},
      author={Picken, Roger},
       title={On two-dimensional holonomy},
        date={2010},
        ISSN={0002-9947},
     journal={Trans. Amer. Math. Soc.},
      volume={362},
      number={11},
       pages={5657\ndash 5695},
}

\bib{MarPickGray}{article}{
      author={Martins, Jo\~ao~Faria},
      author={Picken, Roger},
       title={The fundamental {G}ray 3-groupoid of a smooth manifold and local
  3-dimensional holonomy based on a 2-crossed module},
        date={2011},
        ISSN={0926-2245},
     journal={Differential Geom. Appl.},
      volume={29},
      number={2},
       pages={179\ndash 206},
}

\bib{MarPickSurf}{article}{
      author={Martins, Jo\~ao~Faria},
      author={Picken, Roger},
       title={Surface holonomy for non-abelian 2-bundles via double groupoids},
        date={2011},
        ISSN={0001-8708},
     journal={Adv. Math.},
      volume={226},
      number={4},
       pages={3309\ndash 3366},
}

\bib{Mig1980}{article}{
      author={Migdal, A.~A.},
       title={Properties of the loop average in {QCD}},
        date={1980},
        ISSN={0003-4916},
     journal={Ann. Physics},
      volume={126},
      number={2},
       pages={279\ndash 290},
         url={https://doi-org.libezp.lib.lsu.edu/10.1016/0003-4916(80)90177-3},
      review={\MR{576411}},
}

\bib{Parz2015}{article}{
      author={Parzygnat, Arthur~J.},
       title={Gauge invariant surface holonomy and monopoles},
        date={2015},
        ISSN={1201-561X},
     journal={Theory Appl. Categ.},
      volume={30},
       pages={Paper No. 42, 1319\ndash 1428},
}

\bib{Peif}{article}{
      author={Peiffer, Ren\'ee},
       title={{\"U}ber {I}dentit\"aten zwischen {R}elationen},
        date={1949},
        ISSN={0025-5831},
     journal={Math. Ann.},
      volume={121},
       pages={67\ndash 99},
}

\bib{Pfei2003}{article}{
      author={Pfeiffer, Henryk},
       title={Higher gauge theory and a non-abelian generalization of 2-form
  electrodynamics},
        date={2003},
     journal={Ann. Phys.},
      volume={308},
       pages={447\ndash 477},
}

\bib{SatSchSt}{incollection}{
      author={Sati, Hisham},
      author={Schreiber, Urs},
      author={Stasheff, Jim},
       title={{$L_\infty$}-algebra connections and applications to {S}tring-
  and {C}hern-{S}imons {$n$}-transport},
        date={2009},
   booktitle={Quantum field theory},
   publisher={Birkh\"auser, Basel},
       pages={303\ndash 424},
}

\bib{SchrWalGerb13}{article}{
      author={Schreiber, Urs},
      author={Waldorf, Konrad},
       title={Connections on non-abelian gerbes and their holonomy},
        date={2013},
        ISSN={1201-561X},
     journal={Theory Appl. Categ.},
      volume={28},
       pages={476\ndash 540},
}

\bib{SchreibWalLocal}{article}{
      author={Schreiber, Urs},
      author={Waldorf, Konrad},
       title={Local theory for 2-functors on path 2-groupoids},
        date={2017},
        ISSN={2193-8407},
     journal={J. Homotopy Relat. Struct.},
      volume={12},
      number={3},
       pages={617\ndash 658},
         url={http://dx.doi.org/10.1007/s40062-016-0140-4},
      review={\MR{3691299}},
}

\bib{Sing1995}{incollection}{
      author={Singer, I.~M.},
       title={On the master field in two dimensions},
        date={1995},
   booktitle={Functional analysis on the eve of the 21st century, {V}ol. 1
  ({N}ew {B}runswick, {NJ}, 1993)},
      series={Progr. Math.},
      volume={131},
   publisher={Birkh\"{a}user Boston, Boston, MA},
       pages={263\ndash 281},
}

\bib{SZ}{article}{
      author={Soncini, Emanuele},
      author={Zucchini, Roberto},
       title={A new formulation of higher parallel transport in higher gauge
  theory},
        date={2015},
        ISSN={0393-0440},
     journal={J. Geom. Phys.},
      volume={95},
       pages={28\ndash 73},
}

\bib{Sour1980}{inproceedings}{
      author={Souriau, J.-M.},
       title={Groupes diff\'{e}rentiels},
        date={1980},
   booktitle={Differential geometrical methods in mathematical physics ({P}roc.
  {C}onf., {A}ix-en-{P}rovence/{S}alamanca, 1979)},
      series={Lecture Notes in Math.},
      volume={836},
   publisher={Springer, Berlin-New York},
       pages={91\ndash 128},
}

\bib{stacey2010comparative}{misc}{
      author={Stacey, Andrew},
       title={{C}omparative {S}mootheology},
        date={2010},
}

\bib{WalSurfHol10}{incollection}{
      author={Waldorf, Konrad},
       title={Surface holonomy},
        date={2010},
   booktitle={Handbook of pseudo-{R}iemannian geometry and supersymmetry},
      series={IRMA Lect. Math. Theor. Phys.},
      volume={16},
   publisher={Eur. Math. Soc., Z\"urich},
       pages={653\ndash 682},
}

\bib{Waldorf2016}{misc}{
      author={Waldorf, Konrad},
       title={{ A global perspective to connections on principal 2-bundles}},
         how={\url{https://arxiv.org/abs/1608.00401}},
        date={2017},
}

\bib{Waldorf2017}{misc}{
      author={Waldorf, Konrad},
       title={{ Parallel transport in principal 2-bundles}},
         how={\url{https://arxiv.org/abs/1704.08542}},
        date={2017},
}

\bib{Wang3G}{article}{
      author={Wang, Wei},
       title={On 3-gauge transformations, 3-curvatures, and {\bf
  {g}ray}-categories},
        date={2014},
        ISSN={0022-2488},
     journal={J. Math. Phys.},
      volume={55},
      number={4},
       pages={043506, 32},
}

\bib{Wang2H}{misc}{
      author={Wang, Wei},
       title={{ On the Global 2-Holonomy for a 2-Connection on a 2-Bundle}},
         how={\url{http://arxiv.org/abs/1512.08680}},
        date={2015},
}

\bib{Wang3R}{article}{
      author={Wang, Wei},
       title={On the 3-representations of groups and the 2-categorical traces},
        date={2015},
        ISSN={1201-561X},
     journal={Theory Appl. Categ.},
      volume={30},
       pages={Paper No. 56, 1999\ndash 2047},
}

\end{biblist}
\end{bibdiv}
\bibliographystyle{plain}

\end{document}